\documentclass[twocolumn]{autart}
\usepackage{amssymb,amsfonts,amsmath}
\usepackage{mathtools}
\usepackage{mathrsfs}  
\usepackage{bm,bbm,bbold}
\usepackage{graphicx}
\usepackage{color,ulem}
\usepackage{caption}
\usepackage{subcaption}
\usepackage{epstopdf}                           

\def\dfrac{\displaystyle\frac}

\newtheorem{theorem}{Theorem}[section]
\newtheorem{definition}{Definition}[section]
\newtheorem{lemma}{Lemma}[section]
\newtheorem{remark}{Remark}[section]
\newtheorem{example}{Example}[section]
\renewcommand\thesection{\arabic{section}}
\numberwithin{equation}{section}
\newenvironment{proof}{\begin{pf}}{\rule{1ex}{1ex} \end{pf}}

\newcommand{\R}{\mathbb{R}}
\newcommand{\C}{\mathbb{C}}
\newcommand{\Z}{\mathbb{Z}}

\newcommand{\vsH}{\mathcal{H}}
\newcommand{\vsHza}{\mathring{\mathcal{H}}}

\newcommand{\vsV}{\mathcal{V}}
\newcommand{\vsLt}{L^2}

\newcommand{\vsLiHza}{L^\infty(\R_+,\vsHza)}

\newcommand{\vsVza}{\mathring{\mathcal{V}}}
\newcommand{\DtctCls}{\mathscr{C}_1(h,\Omega)}
\newcommand{\DtctClsA}{\mathscr{C}_2(h,\Omega)}

\newcommand{\ball}{\mathscr{E}_{R}}
\newcommand{\gnrm}{C_{\vg}}

\newcommand{\ve}{\vec e}

\newcommand{\vx}{\vec x}
\newcommand{\vu}{\vec u}
\newcommand{\vz}{\vec z}
\newcommand{\vk}{\vec k}
\renewcommand{\ll}{\ell}
\newcommand{\vf}{\vec f}
\newcommand{\vF}{\vec F}
\newcommand{\vphi}{\vec{\phi}}

\newcommand{\vg}{\vec g}
\newcommand{\vL}{\vec \ell}

\newcommand{\dom}{\mathscr{D}}
\newcommand{\sLCO}{\mathscr{L}}
\newcommand{\oC}{\mathcal C}

\newcommand{\ii}{{i\mkern1mu}}
\newcommand{\ee}{{\mathbb{e}\mkern1mu}}

\begin{document}
\begin{frontmatter}

		\title{Detectability and global observer design for 2D Navier-Stokes equations with uncertain inputs}
		\author[SZ]{Sergiy Zhuk}\ead{sergiy.zhuk@ie.ibm.com},
                \author[SZ]{Mykhaylo Zayats}\ead{mykhaylo.zayats1@ibm.com},
		\author[EF]{Emilia Fridman}\ead{emilia@tauex.tau.ac.il  }
                \address[SZ]{IBM Research - Europe, Dublin, Ireland}
                \address[EF]{ Department of Electrical Engineering-Systems, Tel Aviv University, Israel}

\begin{abstract}
  We present {\it simulation friendly} detectability conditions for 2D Navier-Stokes Equation (NSE) with periodic boundary conditions, and describe a generic class of ``detectable'' observation operators: it includes pointwise evaluation of NSE's solution at interpolation nodes, and spatial average measurements. For ``detectable'' observation operators we design a global infinite-dimensional observer for NSE with uncertain possibly destabilizing inputs: in our numerical experiments we illustrate $H^1$-sensitivity of NSE to small perturbations of initial conditions, yet the observer converges for known and uncertain inputs. 
\end{abstract}
\end{frontmatter}


\section{Introduction}
\label{sec:introduction}
Navier-Stokes Equation (NSE)
\begin{equation}\label{eq:NSEintro}
    \dfrac{d\vu}{dt} + (\vu \cdot \nabla) \vu - \nu\Delta\vu+\nabla p = \vec f
  \end{equation}
  is a basic mathematical model of fluid dynamics: it describes evolution of fluid's velocity vector-field $\vu$ and scalar pressure field $p$ as a function of initial and boundary conditions, input $\vf$ and coefficient $\nu>0$. NSE has applications in biology, weather prediction, energy forecasting to name a few. It also serves as a mathematical model of turbulence: yet an open problem in 3D, in 2D NSE is used to study turbulence and, in particular, to reveal its connection to \textit{deterministic chaos}~\cite{NSE-Turbulence-book2001}. Indeed, for certain $\nu$ and $\vf$ NSE~\eqref{eq:NSEintro} has the unique \textit{attractor}, which is globally exponentially stable, e.g. $\vu=0$ for $\vf=0$ and any $\nu>0$. However, for small enough $\nu>0$ and a destabilising input $\vf$, NSE's attractor could be a multidimensional manifold (e.g.~\cite{NSE-Turbulence-book2001,ilyin2006sharp}), moreover, the attractor could also be chaotic, i.e. initially close-by trajectories might diverge over time on the attractor. For such destabilising inputs classical control problems, e.g. observer design, become challenging, as was noted in~\cite{vazquez2005closed},  especially when attractor's structure depends on an uncertain input $\vf$.

  In this work we generalise the classical notion of detectability for LTI evolution equations~\cite{Bensoussan:Book2007} to periodic 2D NSE, and introduce detectability conditions: we describe a generic class of ``detectable'' observation (output) operators $\oC$ verifying the proposed conditions (our 1st contribution). Intuitively, ``detectable'' output $\oC\vu$ of~\eqref{eq:NSEintro} must provide information about a finite-dimensional subspace where the nonlinear advection term $(\vu \cdot \nabla) \vu$ \textit{is not dominated} by the diffusion term $- \nu\Delta\vu$ (see also the mechanism of energy cascades~\cite{NSE-Turbulence-book2001}). Given such output $\oC\vu(t)$ one can then reconstruct the entire state $\vu(t)$ by designing a minimax filter which injects $\oC\vu$ so that the reconstruction error in the ``unmeasured'' infinite-dimensional orthogonal complement of the range of $\oC$ will decay thanks to the stabilizing effect of $- \nu\Delta\vu$ akin to the case of detectable LTI systems. More specifically, by employing the direct Lyapunov method, we design an \textit{approximate minimax filter} for NSE~\eqref{eq:NSEintro}: we show that $\vu(t)$, the solution of~(\ref{eq:NSEintro}) for an unknown initial condition $\vu(0)=\vu_0$, and for $\vf=\vg+\vec d$, $\vg$ -- a destabilising known input, $\vec d$ -- an uncertain bounded input, belongs to $H^1$-ellipsoid, centered at the state $\vz(t)$ of Luenberger observer, and its radius decays to $0$ as $t\to\infty$ independently of $\vu_0$, provided the uncertain input $\vec d$ ``perturbs'' $\vg$ less frequently in time or $L^2(\Omega)$-norm of the perturbation $\vec d$ (in space) decays as $t\to\infty$, a remarkable result for turbulent systems which are highly sensitive to small perturbations (our 2nd contribution).\\ 
  Infinite-dimensional observers for NSEs, employing backstepping as a design technique, were proposed in~\cite{vazquez2005closed} for Poisseuille flows, a benchmark for turbulence estimation, and in~\cite{he2018observer} for local observer-based stabilization of NSEs. Disturbance estimation by means of backstepping for linear 1D heat equation was proposed in~\cite{feng2017new}. 

  The proposed class of ``detectable'' output operators $\oC$ generalizes the notions of so-called \textit{determining modes}, e.g. spatial averages $\int_{\Omega_j}\vu d\vx$ of the solution $\vu$ over squares $\Omega_j$ covering the domain $\Omega$ (see e.g.~\cite{emilia3,emilia2,wen}), and \textit{determining nodes}, e.g. pointwise evaluation $\vu(\vx_j)$ at finitely many interpolation nodes $\vx_j$ (see e.g.~\cite{azouani2014continuous}), which are well-known in NSE's literature~\cite[p.123,p.131]{NSE-Turbulence-book2001}. Our class contains those cases, and more generally it consists of all closed linear operators $\oC$ verifying certain inequalities: e.g., spatial averages output $\oC\vu = \sum_{j=1}^N \int_{\Omega_j}\vu d\vx \xi_j$, $\xi_j$ -- normalised indicator of $h\times h$-square $\Omega_j$, is detectable if $h^2$, the area of $\Omega_j$ is small enough for $\vu-\oC \vu $ to be ``controlled'' by $h\nabla\vu$ (Poincar\'e inequality). These inequalities do not cover boundary-type observations though (see~\cite{krstic,krstic2}).\\
We stress that spatial averages outputs are important in real-world applications: for example, short-term solar energy forecasting from a sequence of Cloud Optical Depth (COD) images relies upon estimating averaged velocities of clouds from COD images and uses those as measurements in periodic 2D NSE for cloud velocity prediction~\cite{SZTTAA_CDC17}; similarly 2D NSE with averaged velocity measurements obtained from Sea Surface Temperature (SST) satellite images were used in~\cite{herlin2012divergence} to predict short-term SST dynamics.

  In this work we build upon preliminary results presented in our conference papers~\cite{MZCDC21,kang3}: key differences include detectability conditions and more general case of uncertain inputs. Qualitatively similar results were obtained in~\cite{azouani2014continuous} where sufficient conditions for convergence of Luenberger observer were proposed for periodic boundary conditions and known inputs. The authors make use of Brezis inequality to bound the nonlinear advection term in the error equation, and derive conditions for observer's gain and output (e.g. number and size of squares $\Omega_j$ for the case of spatial averages outputs), sufficient for convergence. An experimental assesment of these conditions was recently provided in~\cite{franz2022bleeps}. In contrast, our convergence analysis relies upon a novel one-parametric inequality relating $L^\infty$ and $H^2$-norms of periodic vector-functions, which for certain values of the parameter reduces to Agmon and Brezis inequalities, and S-procedure widely used in Lyapunov stability analysis. As a result, we get \textit{simulation friendly} and less conservative detectability conditions: e.g. for the case of spatial averages outputs our estimate of $h^2$, the area of squares $\Omega_j$, is at least one order of magnitude better for small $\nu>0$, which is of high interest in the case of turbulent flows, and hence for convergence, observer requires averaging the velocity $\vu$ over larger $\Omega_j$ (see Remark~\ref{sec:observer-design-1} below) so less of squares (``sensors'') is needed, and more importantly, the averaging over larger $\Omega_j$ helps reducing the noise hence improving convergence in practice.


\subsection{Mathematical preliminaries}
\label{sec:math-prel-1}
\textbf{Notation. }~Let $\R^n$ denote Euclidian space of dimension $n$ with inner product $\vu\cdot\vec v = \sum_{i=1}^n u_i v_i$, $\R^n_+$ -- non-negative orthant of $\R^n$, $\R^1_+=\R_+$, and for $\vk,\vL\in\R^n$ set $\frac{\vk}{\vL}=(\frac{k_1}{\ell_1}\dots \frac{k_n}{\ell_n})^\top$. Let $\sLCO(H)$ denote the space of all closed linear operators $\oC$ acting in a Hilbert space $H$ with domain $\dom(\oC)\subset H$. The following functional spaces are standard in NSE's theory (~\cite[p.45-p.48]{NSE-Turbulence-book2001}):
\begin{itemize}
\item $L_p^2(\Omega)$ -- space of $\Omega$-periodic functions $u:\Omega\subset\R^2\to\R$ with period $\Omega=(-\frac{\ell_1}2, \frac{\ell_1}2)\times (-\frac{\ell_2}2, \frac{\ell_2}2)$ for some $\ell_{1,2}>0$
  and inner product $(w,v)=\int_\Omega wv dx_1 dx_2$
        \item $L_p^2(\Omega)^2$ -- space of $\Omega$-periodic vector-functions $\vu=\bigl[\begin{smallmatrix}
            u_1 u_2
          \end{smallmatrix}\bigr]$ with inner product $(\vu,\vec\phi) = (u_1,\phi_1)+(u_2,\phi_2)$ and norm $\|\vu\|^2_{\vsLt}=\|u_1\|^2_{\vsLt}+\|u_1\|^2_{\vsLt}$
	\item $H_p^1(\Omega)=\{u\in L_p^2(\Omega): \|\nabla u\|_{\R^2}\in L_p^2(\Omega)\}$,\\ $H_p^1(\Omega)^2=\{\vu=
          \bigl[\begin{smallmatrix}
            u_1 u_2
          \end{smallmatrix}\bigr]\in L_p^2(\Omega)^2: u_{1,2}\in H_p^1(\Omega)\}$ with norm $\|\vu\|^2_{H^1}=\|\vu\|_{\vsLt}^2+\|\nabla\vu\|^2_{\vsLt}$ where $\|\nabla\vu\|^2_{\vsLt}=\int_\Omega \|\nabla \vu\|_{\R^2}^2dx_1dx_2$, and $\|\nabla \vu\|_{\R^2}^2=\|\nabla u_1\|_{\R^2}^2+\|\nabla u_2\|_{\R^2}^2$
	\item $\vsH = \{\vec v \in [L_p^2(\Omega)]^2:\,\nabla\cdot\vec v = 0\}$ -- space of divergence-free vector-functions $\vec v$, $\vsHza = \{\vec v \in \vsH : \quad \int_\Omega \vec v dxdy =0\}$ -- subspace of $\vsH$ of $\vec v$ with zero mean components
        \item $\vsV = \{\vec v\in [H_p^1(\Omega)]^2 :\quad \nabla\cdot\vec v = 0\}$ and $\vsVza = \{\vec v \in \vsV: \quad \int_\Omega \vec v dxdy =0\}$
        \item $L^2(0,T,H)$ -- space of $H$-valued functions $t\mapsto u(t)\in H$ with finite norm $\|u\|^2_{L^2(0,T,H)}=\int_0^T \|u(t)\|_H^2 dt$ for $T\in(0,+\infty)$, e.g. $L^2(0,T,\vsHza)$ -- space of $\vec v(t,x_1,x_2)$ such that $\int_0^T \int_\Omega \|\vec v(t,x_1,x_2)\|_{\R^2}^2dx_1dx_2 <+\infty$ and $\vec v(t,\cdot)$ has zero divergence and zero mean for almost all $t\in (0,T)$
        \item $L^\infty(0,T,H)$ -- space of $H$-valued functions $t\mapsto u(t)\in H$ such that $\|u(t)\|_H \le C<+\infty$ for some $C>0$ and almost all $t\in (0,T)$, $T\in(0,+\infty)$ with finite norm $\|u\|^2_{L^\infty} = \|u\|^2_{L^\infty(0,T,H)}=\min_{C}\{C>0: \|u(t)\|_H \le C\}$
        \end{itemize}
\subsubsection{Bounds for $L^\infty$-norms of periodic vector-functions}
\label{sec:upper-bounds}

The following lemma 
is a key building block used below to define detectability and design observer for NSE. 
\begin{lemma}\label{sec:upper-bounds-linfty}
If $\vu\in H_p^2(\Omega)^2$ and\footnote{This condition is necessary: lemma does not hold for $u\equiv \text{const}$ and small enough $\ell_{1,2}$} $\int_\Omega \vu dx_1dx_2=0$ then for any $\gamma>0$ it holds:
\begin{equation}
  \label{eq:AgmoN}
  \begin{split}
  \|\vu\|_{L^\infty}&\le \frac{\log^{\frac12}\bigl(1+\frac{4\pi^2 \gamma^2}{\ll_1\ll_2}\bigr)\|\vu\|_{H^1}}{\sqrt{2\pi}} + \frac{\|\vL\|_{\R^2} \|\Delta \vu\|_{L^2}}{\gamma\sqrt{32\pi^3}}
\end{split}
\end{equation}
For $\gamma= \|\vu\|_{H^1}^{-\frac12}\|\Delta \vu\|_{L^2}^{\frac12}$ \eqref{eq:AgmoN} gives 2D Agmon inequality~\cite[p.100]{NSE-Turbulence-book2001}:
\begin{equation}
  \label{eq:Agmon}
  \begin{split}
  \|\vu\|_{L^\infty}&\le (\sqrt{\frac{2\pi}{\ll_1\ll_2}}+\frac{\|\vL\|_{\R^2}}{\sqrt{32\pi^3}})\|\vu\|_{H^1}^{\frac12}\|\Delta \vu\|_{L^2}^{\frac12}
\end{split}
\end{equation}
For $\gamma= \|\vu\|_{H^1}^{-1}\|\Delta \vu\|_{L^2}$ \eqref{eq:AgmoN} gives 2D Brezis inequality~\cite{brezis1979nonlinear}:
\begin{equation}
  \label{eq:Brezis}
  \|\vu\|_{L^\infty}\le \bigl(\frac{\|\vL\|_{\R^2}}{\sqrt{32\pi^3}}+\frac{\log^{\frac12}(1+\frac{4\pi^2 \|\Delta \vu\|_{L^2}^2}{\ll_1\ll_2 \|\vu\|_{H^1}^2})}{\sqrt{2\pi}}\bigr)\|\vu\|_{H^1}
\end{equation}
\end{lemma}
The proof is given in the appendix. 
\subsubsection{Navier-Stokes equation: weak  formulation and well-posedness in 2D}
\label{sec:navi-stok-equat}
 The classical NSE in 2D is a system of two PDEs defining dynamics of the scalar pressure field $p(x,y)$ and the viscous fluid velocity vector-field $\vu(t,x,y)=[\begin{smallmatrix} u_1(t,x,y),u_2(t,x,y)
   \end{smallmatrix}]$ which depends on the initial condition $\vu(0)=\vu_0\in\vsH$, input (e.g. forcing) $\vf=[f_1,f_2]$ and Boundary Conditions (BC), e.g. periodic BC $u_{1,2}(t,x+\ell_1,y) = u_{1,2}(t,x,y)$, $u_{1,2}(t,x,y+\ell_2) = u_{1,2}(t,x,y)$. In the vector form it reads as follows:
\begin{equation}
  \label{eq:NSE-vector}
  \begin{split}
    \dfrac{d\vu}{dt} &+ (\vu \cdot \nabla) \vu - \nu\Delta\vu+\nabla p = \vf,\quad \nabla\cdot\vu = 0
  \end{split}
\end{equation}
To eliminate pressure $p$ and obtain an evolution equation just for $\vu$ it is common to use Leray projection~\cite[p.38]{NSE-Turbulence-book2001}: every vector-field $\vu$ in $\R^2$ admits Helmholtz-Leray decomposition, $\vu = \nabla p + \vec v$ with $\nabla\cdot\vec v = 0$ which in turn defines Leray projector, $P_l(\vu)=\vec v$ -- an orthogonal projector onto $\vsHza$ (e.g.~\cite[p.36]{NSE-Turbulence-book2001}). 
Multiplying~\eqref{eq:NSE-vector} by a test function $\vec\phi\in\vsVza$ (the projection step), and integrating by parts in $\Omega$ allows one to obtain Leray's weak formulation of NSE in 2D:
	\begin{equation}
	\label{eq:NSE-var}
	\dfrac{d}{dt}(\vu,\vec \phi) + b(\vu,\vu,\vec\phi) + \nu((\vec u,\vec \phi)) = (\vec f,\vec \phi), \quad\forall\vec\phi\in \vsVza
      \end{equation}
      with initial condition $(\vu(0),\vec\phi)=(\vu_0,\vec\phi)$. Here
\begin{align*}
 b(\vec u,\vec w,\vec\phi) &= (\vec u\cdot \nabla w_1,\phi_1) + (\vec u\cdot \nabla w_2,\phi_2),\\
                             ((\vec u,\vec \phi)) &= (\nabla u_1,\nabla \phi_1)+(\nabla u_2,\nabla \phi_2)
\end{align*}
In what follows we will be using some properties of the trilinear form $b$ and Stokes operator  $\vu\mapsto A\vu=-P_l\Delta \vu$, a self-adjoint positive operator with compact inverse, which coincides with $\Delta \vec u$ for periodic BC (see~\cite[p.52]{NSE-Turbulence-book2001}): for $\vu\in \dom(A)$ and $\vphi\in \vsVza$
\begin{align}
  (A\vec u, \vphi) &= ((\vec u, \vphi)),\,(A\vec u,\vec u) = ((\vec u,\vec u))\ge \lambda_1(\vu,\vu)  \label{eq:Poincare}\\
  (A\vec u, A\vec u)  &=(A(A^\frac12)\vec u,(A^\frac12)\vec u)\ge \lambda_1 (A \vec u, \vec u)\label{eq:gradA-bnd}\\
  \lambda_1&=4\pi^2/\max\{\ell_1,\ell_2\}^2\\
 b(\vec u,\vec v,\vec \phi) &= - b(\vec u,\vec \phi,\vec v)\\
 b(\vec u,\vec v,\vec v) &=0,\quad  b(\vec v,\vec v,A\vec v) = 0 \label{eq:b:ort}
\end{align}

Next Lemma collects results from~\cite[p.58, Th.7.4, p.99, f.(A.42), p.102, f.(A.66)-(A.67)]{NSE-Turbulence-book2001} on existence, uniqueness, regularity and input-to-state stability of NSE's weak (strong) solution $\vu$, and bounds for $A\vu$. Classical smoothness of $\vu$ requires further constraining of $\vf$, $\vu_0$~\cite[p.59]{NSE-Turbulence-book2001}. 
\begin{lemma}
	Let $\vu_0\in \vsHza$ and $\vf\in L^2(0,\bar T,\vsHza)$. Then, on $[0,\bar T]$ there exist the unique weak solution $\vu\in C(0,\bar T,\vsHza)$ of NSE~\eqref{eq:NSE-var}, and the components of $\vu=\begin{smallmatrix}[u_1 & u_2]\end{smallmatrix}$ verify: $u_i,(u_i)_{x,y}\in L^2(\Omega\times(0,\bar T))$. If $\vu_0\in\vsVza$ then the weak solution coincides with the strong solution of~\eqref{eq:NSE-var}, and \[
	\dfrac{du_i}{dt}, (u_i)_{x},(u_i)_{y}, (u_i)_{x y}\in L^2(\Omega\times(0,\bar T))
      \] i.e. $\vu\in C(0,\bar T,\vsVza) \cap L^2(0,\bar T,\mathscr{D}(A))$. If in addition $\vec f\in \vsLiHza$ then $(\,\|\vec f\|_{L^\infty} = \|\vec f\|_{\vsLiHza}$ for short$)$:
	\begin{align}
          &\|\vec u(\cdot,t)\|^2_{\vsLt} \le \frac{\|\vec f\|^2_{L^\infty}}{(\nu\lambda_1)^2} + e^{(-\lambda_1 \nu)(t-s)}\|\vec u(\cdot,s)\|^2_{\vsLt}	\label{eq:L2norm-decay}\\
	&\|\nabla \vec u(\cdot,t)\|^2_{\vsLt} \le \frac{\|\vec f\|^2_{L^\infty}}{\nu^2\lambda_1} + e^{(-\lambda_1 \nu)(t-s)}\|\nabla \vec u(\cdot,s)\|^2_{\vsLt}          	\label{eq:grad-norm-decay}\\
          &\frac1{\bar T}\int_t^{\bar T+t}\|A\vu\|^2_{\vsLt}ds  \le \theta_{t,\bar T}:=\frac{2\|\vf\|^2_{L^\infty}}{\bar T\nu^3\lambda_1} + \label{eq:Au-L2avrg-bnd}\\
          &\quad + \frac1{\bar T}\int_t^{\bar T+t}\frac{\|\vf(s)\|^2_{\vsLt}}{\nu^2}ds + \frac{2 e^{(-\lambda_1 \nu)t}\|\nabla\vec u_0\|^2_{\vsLt}}{\bar T\nu}\label{eq:Au-ThetatT}
        \end{align}
\end{lemma}


\section{Problem statement}
\label{sec:problem-statement}
Generalizing the classical definition of detectability for LTI systems~\cite{Bensoussan:Book2007} we say that NSE is detectable w.r.t. output operator $\oC$ if the distance (in space) between any two solutions corresponding to different initial conditions and inputs converges to zero over time, provided so does the distance between their respective outputs and inputs:
\begin{definition}\label{def:detectability}
  For $\vf,\vF\in \vsLiHza$ let $\vu$ and $\vz$ solve
  \begin{align}
    \dfrac{d}{dt}(\vec u,\vec \phi) &+ b(\vu,\vu,\vec\phi) + \nu((\vu,\vec \phi)) = (\vec f,\vec \phi)\,,\vu(0)\in\vsVza\label{eq:NSEu}\\
\dfrac{d}{dt}(\vec z,\vec \phi) &+ b(\vec z,\vec z,\vec\phi) + \nu((\vec z,\vec \phi)) = (\vec F,\vec \phi)\,, \vz(0)\in\vsVza \label{eq:NSEz}
  \end{align}
NSE is called detectable in $\vsVza$ w.r.t. an linear output operator $\oC\in\mathscr{L}(\vsH)$ if the following conditions 
  \begin{align*}
    &\mathrm{A)}~\exists T>0:\limsup_{t\to\infty}\frac1T\int_t^{t+T}\|\oC(\vu(\cdot,s)-\vz(\cdot,s))\|_{\vsLt}^2 ds=0\\
    &\mathrm{B)}~\exists T>0:\limsup_{t\to\infty}\frac1T\int_t^{t+T} \|\vec f(\cdot,s)-\vec F(\cdot,s)\|_{\vsLt}^2 ds =0
  \end{align*}
imply state convergence in $\vsVza$: $\|\vu(\cdot,t)-\vz(\cdot,t)\|_{\vsVza}\underset{t\to\infty}{\to}0$.
\end{definition}
In practice, typical output operators $\oC$ are either pointwise evaluations (see Example~\ref{exmpl:detectability-nse-2} below), or of projection type, i.e. $\oC\vu$ is a projection of $\vu$ onto a subspace (see Example~\ref{exmpl:detectability-nse-1} below). The following definition introduces two classes of output operators $\oC$ which generalize projection type and pointwise evaluation type observations:
\begin{definition}\label{def:detectability-class}
 Take $\oC\in\sLCO(\vsLt)$ such that $\vsV\subseteq\dom(\oC)$ and fix $h>0$. $\oC$ is said to be of class $\DtctCls$ or of class $\DtctClsA$ if~\eqref{eq:CPoin} or ~\eqref{eq:CPoin2} respectively holds for some $C_\Omega>0$
  \begin{align}
    \forall\vu\in\vsVza:\quad \|\vu - \oC\vu\|_{\vsLt}^2&\le h^2 C_\Omega \|\nabla\vu\|_{\vsLt}^2 \label{eq:CPoin}\\
    \forall\vu\in \dom(A):\quad \|\vu - \oC\vu\|_{\vsLt}^2&\le h^2 C_\Omega \|\Delta\vu\|_{\vsLt}^2 \label{eq:CPoin2}
  \end{align}
\end{definition}
\begin{example}\label{exmpl:detectability-nse-1}
  The class $\DtctCls$ generalizes the notion of so-called \textit{determining modes} for NSE (e.g.~\cite[p.123]{NSE-Turbulence-book2001}). The case of determining modes corresponds to $\oC\in\DtctCls$ of finite rank, $\oC\vu$ represents the projection of $\vu$ onto a certain finite-dimensional subspace of $\vsH$, and in this case $\oC$ is bounded in $\vsH$. As an example, consider the case of spatial averages over a partition of $\Omega$ (e.g.~\cite{emilia2,azouani2014continuous}):  Let $\Omega=\cup_{j=1}^N \Omega_j$, $\Omega_j\cap\Omega_i=\varnothing$ for $i\ne j$, $\Omega_j$ -- a rectangle with sides of length $h_x^j$ and $h_y^j$ respectively, and let $\xi_j$ denote the indicator function of $\Omega_j$. 
  Let \[
    \oC\vu(x,y) = \sum_{j=1}^N\bigl[
\begin{smallmatrix}
  (u_1,\xi_j)\\(u_2,\xi_j)
\end{smallmatrix}\bigr]\frac{\xi_j(x,y)}{h_x^j h_y^j}, \quad \xi_j(x,y)=1, (x,y)\in\Omega_j
\] Then~\eqref{eq:CPoin} holds for $h^2=\max_j(\max\{h_x^j,h_y^j\})^2$ and $C_\Omega=(4\pi^2)^{-1}$ by Poincar\'e inequality~\eqref{eq:Poincare}.
\end{example}
\begin{example}\label{exmpl:detectability-nse-2}
  The class $\DtctClsA$ generalizes notion of so-called \textit{determining nodes} for NSE (see e.g.~\cite[p.131]{NSE-Turbulence-book2001}). The case of determining nodes corresponds to $\oC\in\DtctClsA$ of the form  $\oC\vu=\sum_{j=1}^n \vu(t,\vx_j)\vphi_j$ for some grid nodes $\vx_j\in\Omega$ and functions $\vphi_j\in H^1_p(\Omega)^2$. Clearly, $\oC$ is not bounded in $\vsH$ as it is defined for $\vu\in H^2_p(\Omega)^2\cap\vsH$, but it is a closed linear operator such that $\dom(A)\subset\dom(C)$. As an example, consider the case of pointwise evaluations over nodes $\vx_j$ of a finite element grid, here $\vphi_j$ are 2D hat basis functions. In this case $\oC$ verifies~\eqref{eq:CPoin2} with constants $h$ and $C_\Omega$ which can be found in~\cite[Th 12.3.4]{Hutson1980}.
\end{example}
With the above two definitions in mind we are ready to state \emph{goals of the paper:}
\begin{itemize}
\item \textbf{detectability:} find sufficient conditions for detectability (as per Definition~\ref{def:detectability}) of NSE w.r.t. any output operator $\oC\in\DtctCls$ or $\oC\in\DtctClsA$
\item \textbf{observer design:} given output $\oC\vu$ of~\eqref{eq:NSEu} to construct an output feedback control $\vec F$ such that conditions $A)$ and $B)$ hold and so, by detectability, $\vz$ converges to $\vu$ in $\vsVza$ independently of $\vu(0)$ and $\nabla\vu(0)$
\end{itemize}


\section{Main results}
\label{sec:main-results}
\subsection{Sufficient conditions for detectability}
\label{sec:detectability-nse}
Recall the classes $\DtctCls$ and $\DtctClsA$ of observation operators $\oC$ introduced above in Def.~\ref{def:detectability-class}. The following proposition demonstrates \emph{existence} of a constant $h>0$ for each of $\DtctCls$ and $\DtctClsA$ such that NSE is detectable (as per Definition~\ref{def:detectability}). However, as is, this proposition is of rather theoretical value as to compute $h>0$ for numerical simulations one needs to know certain interpolation constants which are either unknown or very concervative. The question of how to use this result for estimating $h$ in numerical simulations will be addressed below.
\begin{prop}\label{prop:detectability-nse-1}
  Recall Def.~\ref{def:detectability}: let $\vz,\vu$ solve~\eqref{eq:NSEu}-\eqref{eq:NSEz} for $\vf,\vF\in \vsLiHza$. Recall from~\cite[p.100,(A.47)]{NSE-Turbulence-book2001} an interpolation inequality for $\nabla\vu $, and define constants $C_{1,2}$:
\begin{align}
  &\forall\vu\in\dom(A):\|\nabla\vu\|^2_{\vsLt}\le C_{\nabla} \|\vu\|_{\vsLt} \|A\vu\|_{\vsLt}   \label{eq:InterpIneq}\\
  &C_1^2=(1+\lambda_1^{-1})/(2\pi), \quad C_2^2 = \|\vL\|^2_{\R^2}/(32\pi^3)\label{eq:C1C2}
\end{align}
Take $\oC\in\DtctCls$. Then NSE is detectable in $\vsVza$ if for some $1<\kappa<<2$ and $\Gamma>0$
  \begin{equation}
    \label{eq:h-bnd-detect}
    h< \frac{\lambda_1^{\frac12}\nu\bigl(\frac{3\nu}{4} -\frac{C_2}{\lambda_1\Gamma}\bigr)}{C_1 C_\nabla C^{\frac12}_\Omega \kappa^{\frac12}\|\vf\|_{\vsLiHza}\log^{\frac12}\bigl(1+\frac{4\pi^2 \kappa \|\vf\|_{\vsLt}^2\Gamma^2 }{\nu^2\ll_1\ll_2}\bigr)}
  \end{equation}
  moreover, if $\oC\in\DtctClsA$ then NSE is detectable in $\vsVza$ if $h$ verifies~\eqref{eq:h-bnd-detect} without $\lambda_1^{\frac12}$ in the numerator.
\end{prop}
The proof is provided in the appendix.


\subsection{Observer design}
\label{sec:observer-design}
As noted above, 
  condition~\eqref{eq:h-bnd-detect} is hard to use in simulations. Below we build on~\eqref{eq:h-bnd-detect} and propose ``simulation friendly'' conditions for $h$ (see $C1,C2$ of Theorem~\ref{t:1}): it is demonstrated that plugging Luenberger-type output feedback $\vF =\vg+L\oC(\vu-\vz)$ into~\eqref{eq:NSEz} ensures output and input convergence (conditions $A)$ and $B)$ of Def.~\ref{def:detectability}), and as a result implies state convergence in $\vsVza$: $\|\vu(\cdot,t)-\vz(\cdot,t)\|_{\vsVza}\to0$.\\
\begin{lemma}[Wellposedness]\label{l:observer-existence}
  Let $\vu$ solve~\eqref{eq:NSEu}. For any $\vg\in\vsLiHza$ there is the unique $\vz\in L^\infty(\R_+,\vsVza)\cap L^2(t_0,t_1,\dom(A))$, $0\le t_0<t_1<+\infty$ such that
  \begin{align}
    &\dfrac{d}{dt}(\vec z,\vec \phi) + b(\vec z,\vec z,\vec\phi) + \nu((\vec z,\vec \phi)) = (\vec F,\vec \phi)\,, \phi\in\vsVza\label{eq:NSE-LO}\\
&\vec F = \vg+L\oC(\vu-\vz) \label{eq:observerF}
  \end{align}
  provided $\oC\in\DtctCls$ and $L h^2 C_\Omega\le 2\nu$, or $\oC\in\DtctClsA$ and $L hC_\Omega^\frac12\le \nu$. 
\end{lemma}
The proof is given in the appendix.
\begin{theorem}\label{t:1}
  Assume that (i) $\vu$ solves~\eqref{eq:NSEu} for unknown $\vu(0)\in\vsVza$ and $\vf$, and (ii) $\vf=\vg+\vec d$ for a known $\vg\in\vsLiHza$, $\gnrm=\|\vg\|_{\vsLiHza}$ and an uncertain input $\vec d$ from a ball $\ball=\{\vec d:\|\vec d\|_{\vsLiHza} \le R\}$ of radius $R>0$. Let $\vz$ solve~\eqref{eq:NSE-LO} for $\vz(0)=0$, $L>0$ and $\vg$ as in (ii) above. For $\varepsilon>0$ define $\kappa=1+\varepsilon$ and functions of a parameter $\Gamma$:
  \begin{align}
\hat L_\nabla (\Gamma)& = 2(\nu -\frac{\|\vL\|_{\R^2}}{\sqrt{32\pi^3}\lambda_1\Gamma})/C_\Omega, \quad \vL=
   \bigl[\begin{smallmatrix}
     \ell_1\\\ell_2
   \end{smallmatrix}\bigr]\label{eq:hatL}\\
    \hat L_\Delta(\Gamma)&= \frac{(\kappa+\kappa\lambda_1^{-1})^\frac12 \gnrm}{\nu\sqrt{2\pi}}  \log^{\frac12}\bigl(1+\frac{4\pi^2 \kappa \,\gnrm^2\, \Gamma^2 }{\nu^2\ll_1\ll_2}\bigr) \label{eq:hatL1}\\
                   \Theta(\Gamma)&= \hat L_\nabla(\Gamma) / (2 \hat L_\Delta(\Gamma)) \label{eq:Theta}
  \end{align}
and let $\Gamma_{\mathrm{max}}>0$ maximize $\Theta(\Gamma)$. Finally assume that 
  \begin{equation}
    \label{eq:fgdecay}
    \exists T>0:\Sigma_T(s) = \sup_{t\ge s}\frac1{T}\int_t^{t+T}\|\vec d(\tau)\|^2_{\vsLt} d\tau\xrightarrow{s\to\infty}0
  \end{equation}
  and take minimal $t^\star>0$ and $T_1\ge T$ verifying
  \begin{equation}
    \label{eq:tstarTeps}
    \begin{split}
    &\delta(t^\star)+\frac{2(\gnrm+R)^2}{\gnrm^2 T_1\nu\lambda_1} + \frac{2\nu\mathbb{e}^{-\lambda_1 \nu t^\star}\|\nabla\vec u_0\|^2_{\vsLt}}{T_1 \gnrm^2}
    \le \varepsilon\\
    &\delta(t^\star) = \gnrm^{-2}\Sigma_{T_1}(t^\star) + \gnrm^{-1}\Sigma_{T_1}^{\frac12}(t^\star)
  \end{split}
\end{equation}
Then, for $t\ge s\ge t^\star$, $V(t)=\|\nabla(\vu(t)-\vz(t))\|^2_{\vsLt}$ verifies
\begin{equation}
  \label{eq:ineq-V1}
  \begin{split}
    V(t) &\le e^{\omega LT_1} V(s) e^{(Q(L)-\omega L)(t-s) + \omega L T_1 \lfloor\frac{t-s}{T_1}\rfloor}\\
    &+\frac{T_1 e^{T_1 (\omega L - Q(L))}}{C_\Omega(1-\beta) \hat L_\nabla(\Gamma_{\mathrm{max}})} \Sigma_{T_1}(s) \\
    Q(L)&= -L\omega  + \hat L_\Delta(\Gamma_{\mathrm{max}})<0,
\end{split}
\end{equation}
if $\oC$, $L$ in~\eqref{eq:NSE-LO}, and $\beta,\omega$ in~\eqref{eq:ineq-V1} satisfy either $C1$ or $C2$: 
 \begin{itemize}
 \item [$C1)$]$\oC\in\DtctCls$ verifies~\eqref{eq:CPoin} with $h^2<\beta \Theta(\Gamma_{\mathrm{max}})$ for some $0<\beta<1$, and $L=\beta \hat L_\nabla(\Gamma_{\mathrm{max}})/h^2$, $\omega= 1/2$
 \item [$C2)$] $\oC\in\DtctClsA$ verifies~\eqref{eq:CPoin2} with $h< \beta C_\Omega^{\frac12}\theta^{-1}\Theta(\Gamma_{\mathrm{max}})$ for some $0<\beta<1$, $\theta>1$, and $L=\theta\hat L_\Delta(\Gamma_{\mathrm{max}})$, $\omega=1$
 \end{itemize}
\end{theorem}
The proof is provided in the appendix.
\begin{remark}
  Condition~\eqref{eq:fgdecay} is verified if either (i) the average of $L^2$-energy of the uncertain input $\vec d$ over time window $(t,t+T)$ decays to $0$, or (ii) the measure of time instants $s$ within a ``window'' $(t,t+T)$, where $\vec d(s,\cdot)$ is ``active'', decays to $0$ as the window $(t,t+T)$ slides to infinity ($t\to\infty$). In other words, the case (ii) does not require $\|\vec d(s,\cdot)\|_{\vsLt}\to 0$ but requires that the uncertain input $d$ ``perturbs'' $\vg$ less and less frequently asymptotically. \\
  The rate of decay of $\Sigma_T(s)$ to $0$ determines the rate of decay of $V$: indeed, \eqref{eq:ineq-V1} implies that \[
    \sup_{t>\bar t(s)>s} V(t)\le 2 \frac{T_1 e^{T_1 (\omega L - Q(L))}}{C_\Omega(1-\beta) \hat L_\nabla(\Gamma_{\mathrm{max}})} \Sigma_{T_1}(s)
    \] provided $\bar t(s)$ is large enough for the impact of the 1st term in r.h.s. of~\eqref{eq:ineq-V1} be negligible compared to the 2nd term. If $\Sigma_T(s)=0$ for $s>s^\star$ then $V(t)\to 0$ exponentially after $t=\max\{s^\star,t^\star\}$. In fact, $V\le V_1=\|\vu-\vz\|_{H^1(\Omega)}$ and by~\eqref{eq:Poincare} $V_1\le (1+\lambda_1^{-1}) V$ hence theorem proposes an approximation of minimax filter in the following sence: by \eqref{eq:ineq-V1} $\vu(t)$ belongs to $H^1$-ellipsoid, centered at $\vz(t)$, and its radius is given by r.h.s. of~\eqref{eq:ineq-V1}. 
\end{remark}
\begin{remark}\label{sec:observer-design-0}
   Parameter $0<\beta<1$ in $C1$ and $C2$ allows to balance ``the number of sensors'' (e.g. for outputs in the form of spatial averages  larger $h^2$ means less data as one needs less of $h\times h$-squares to cover $\Omega$), and the impact of $\Sigma_T(s)$ onto the decay of $V$: $\beta$ close to $1$ allows to take $h^2$ close to its maximal value of $\Theta(\Gamma_{\mathrm{max}})=\max_{\Gamma>0} \Theta(\Gamma)$ at the price of amplifying the impact of the term with $\Sigma_T(s)$ in~\eqref{eq:ineq-V1} proportionally to $1/(1-\beta)$. If $\Sigma_T(s)=0$ for $s>s^\star$ one can set $\beta=1$. Parameter $\theta>1$ is relevant only for outputs of $\DtctClsA$-class, it is independent of $\beta$ and allows one to balance convergence speed and amount of observations: larger $\theta$ will imply faster convergence rate at the price of decreasing $h$ proportionally to $1/\theta$. Interestingly, $L$ increases with the decrease of $h$ as per $C1$ but is independent of $h$ if the case $C2$.
\end{remark}
\begin{remark}\label{sec:observer-design-1}
  Let us compare our upper bound for $h^2$ obtained above for $\oC\in\DtctCls$, namely $h^2<\beta\Theta(\Gamma_{\mathrm{max}})$ to the state-of-the-art result of~\cite[Prop.2]{azouani2014continuous}:
  \begin{equation}
    \label{eq:Titih}
    h^2\le \nu \mathcal{T}(\nu)=C_\Omega^{-1} \bigl(3\nu\lambda_1 (2c\log(2) c^{\frac32}+8c\log(1+G))G\bigr)^{-1},
  \end{equation}
  where $G=\frac{\|\vf\|_{\vsLiHza}}{\lambda_1\nu^2}$ and the constant $c$ comes from Brezis inequality $\|\vu\|_{L^\infty}\le c \|\nabla\vu\|_{\vsLt}(1+\log\frac{\|A\vu\|_{\vsLt}^2}{\lambda_1\|\nabla\vu\|_{\vsLt}^2})^{\frac12}$ for $\vu\in\dom(A)$~\cite[p.100, (A.50)]{NSE-Turbulence-book2001}. Clearly, \[
    c=\sup_{\vu\in\dom(A):\|\nabla\vu\|_{\vsLt}>0}\|\vu\|_{L^\infty} \|\nabla\vu\|^{-1}_{\vsLt}(1+\log\frac{\|A\vu\|_{\vsLt}^2}{\lambda_1\|\nabla\vu\|_{\vsLt}^2})^{-\frac12}
  \] and so $c\ge c_1$ for $c_1=\|\vu_1\|_{L^\infty} \|\nabla\vu_1\|^{-1}_{\vsLt}(1+\log\frac{\|A\vu_1\|_{\vsLt}^2}{\lambda_1\|\nabla\vu_1\|_{\vsLt}^2})^{-\frac12}$ and $\vu_1\in\dom(A)$ with components $u_1=-\cos(2\pi x)\sin(2\pi y)$ and $u_2=\cos(2\pi y)\sin(2\pi x)$. It is easy to compute that $c_1=(2\pi)^{-1}$. Since the r.h.s. of~\eqref{eq:Titih}, $\mathcal{T}(\nu)$ increases if we substitute $c$ with $c_1$, and so the upper bound on $h^2$ improves, below we compare $\Theta(\Gamma_{\mathrm{max}})$ vs. $\mathcal{T}(\nu)$ with $c=c_1=(2\pi)^{-1}$ over the interval $\nu\in[10^{-6},10^{-1}]$. To match the setting of~\cite{azouani2014continuous} we assume that $\ell_{1,2}=1$, $\|\vf(s)-\vg(s)\|_{\R^2}=0$ for $s>s^\star$ so we can use $\beta=1$, we also take $\|\vf\|_{\vsLiHza} = C_{\vf}=1$ and $\kappa=2$. We get that $\log_{10}\frac{\Theta(\Gamma_{\mathrm{max}})}{\mathcal{T}(\nu)}=1.33$ for $\nu=10^{-6}$ and $0.58$ for $\nu=10^{-1}$ and LogLog-plot of $\Theta$ and $\mathcal{T}$ over $\nu\in[10^{-6},10^{-1}]$ is given in Fig.~\ref{fig:comp}. Obviously, for small $\nu$ our upper bound is at least one order of magnitude better.
 \begin{figure}[h]
\begin{center}
  \includegraphics[width=\columnwidth]{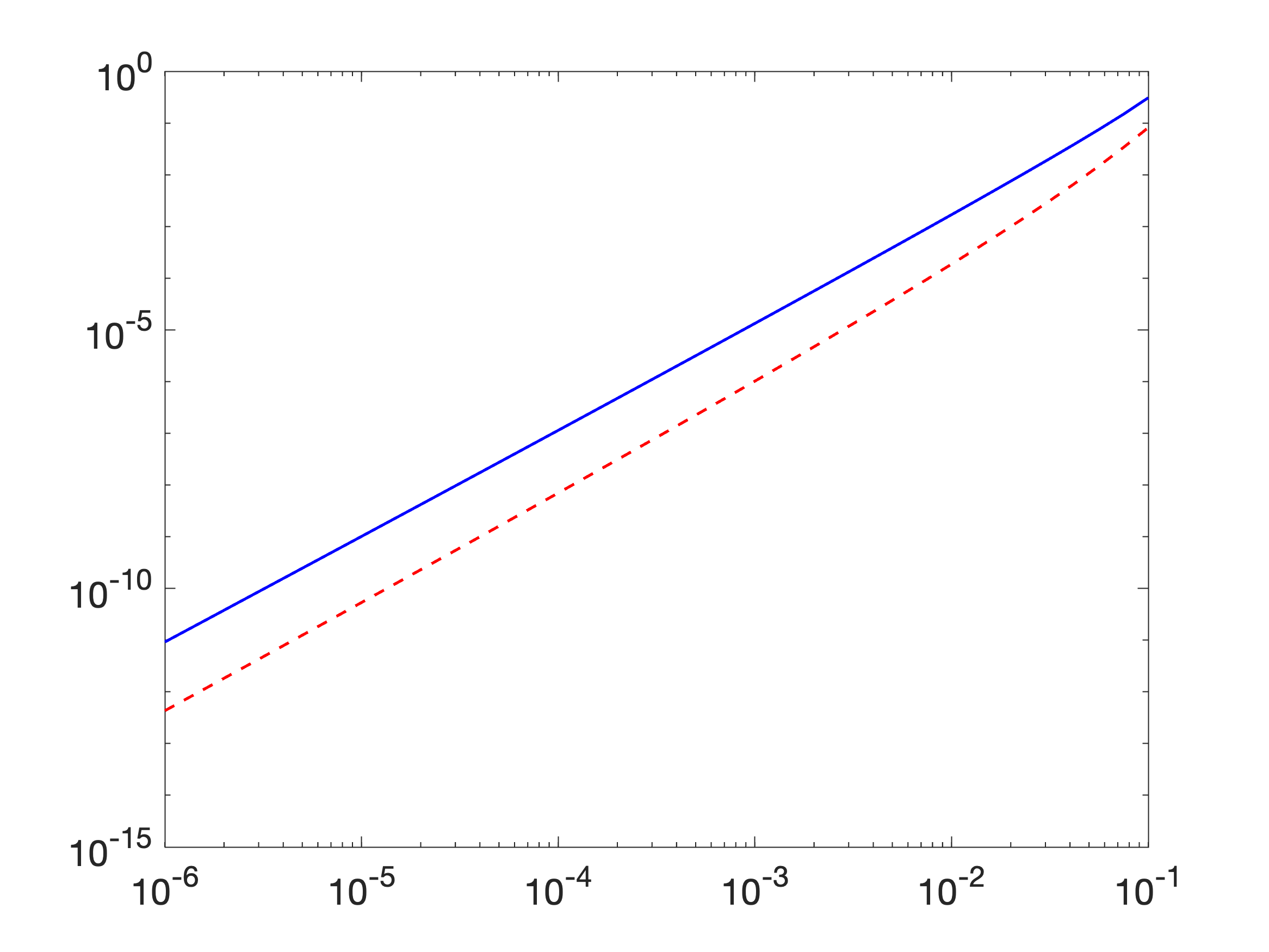}
  \caption{LogLog-plot of upper-bounds for $h^2$: $\Theta(\Gamma_{\mathrm{max}})$ (blue) vs. $\mathcal{T}(\nu)$ (red) with $c=c_1=(2\pi)^{-1}$ for $\nu\in[10^{-6},10^{-1}]$.}
	\label{fig:comp}
	\end{center}
\end{figure}
\end{remark}


\section{Experiments}
\label{sec:experiments}

\begin{figure}
            \centering
\includegraphics[width=\columnwidth]{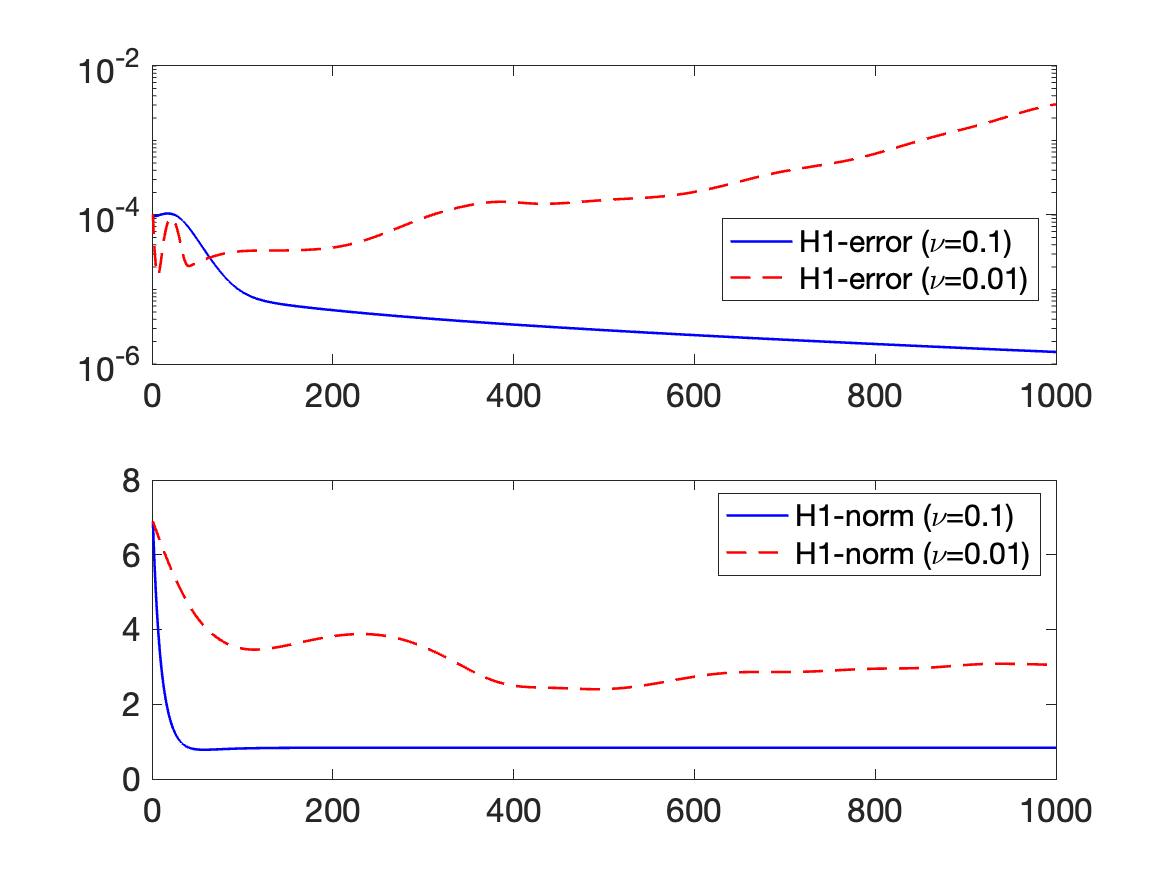}
\caption{Top panel: $H^1$-relative difference of initially close-by trajectories for different $\nu$ over time. Bottom panel: $H^1$-norm dynamics of solution with the same initial condition and input for different $\nu$ over time.}
\label{fig:TurbulenceA}
\end{figure}
\begin{figure}
\captionsetup[subfigure]{aboveskip=-10pt,belowskip=-3pt}
\begin{subfigure}{\columnwidth}
  \centering
  \includegraphics[width=\columnwidth]{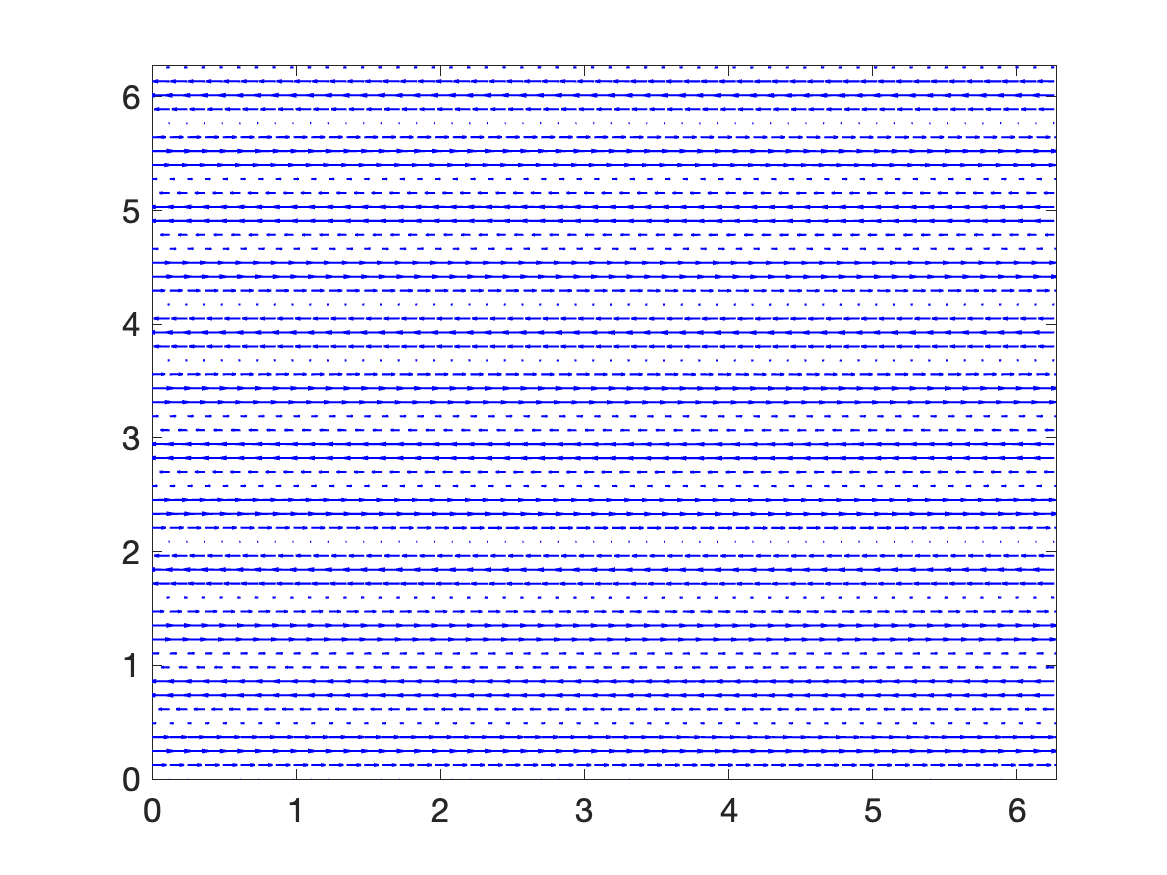}
  \caption{Example of non-turbulent velocity for $\nu=0.1$ (at final time) }
  \label{fig:nu1}
\end{subfigure}
\captionsetup[subfigure]{aboveskip=-10pt,belowskip=-10pt}
\begin{subfigure}{\columnwidth}
  \centering
  \includegraphics[width=\columnwidth]{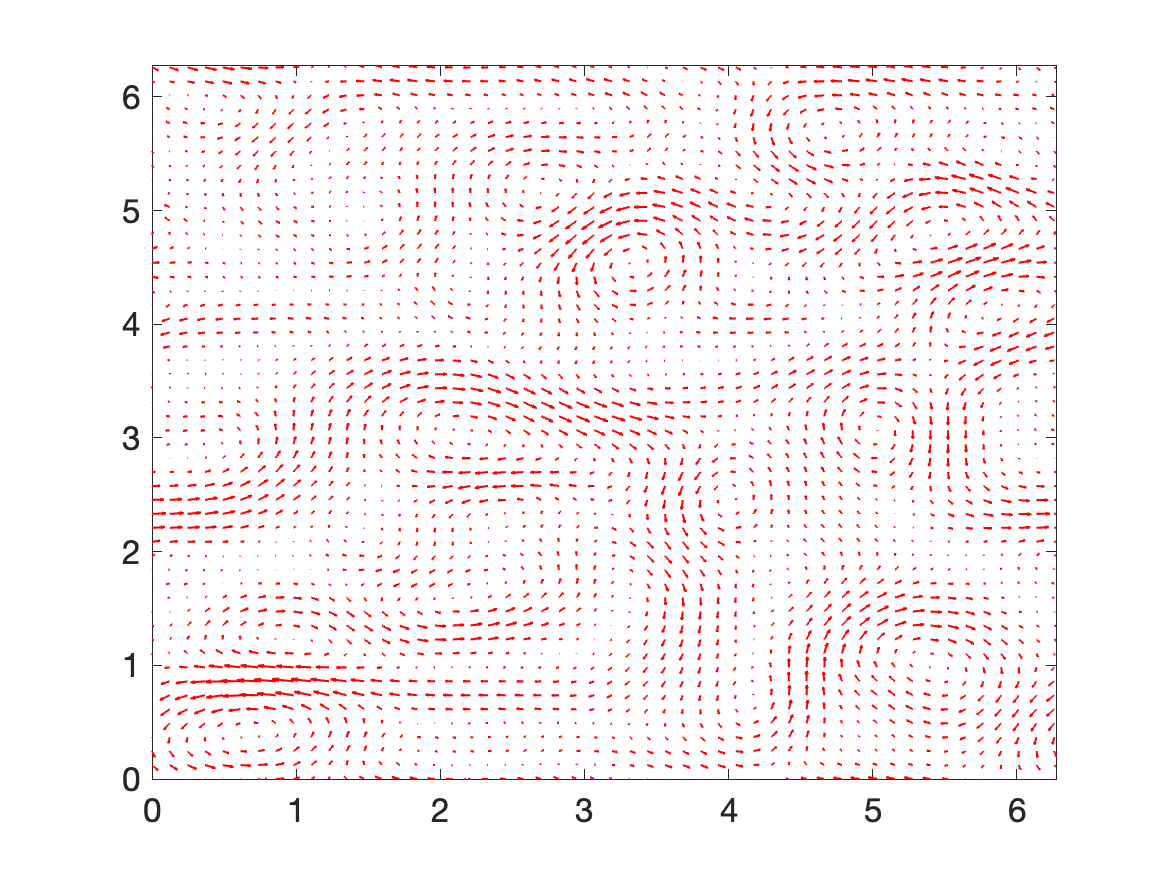}
  \caption{Example of turbulent velocity for $\nu=0.01$ (at final time) }
  \label{fig:nu01}
\end{subfigure}
\caption{Turbulence: velocity snapshots}
\label{fig:Turbulence}
\end{figure}
\begin{figure}
 \centering 
 \includegraphics[width=1.0\columnwidth]{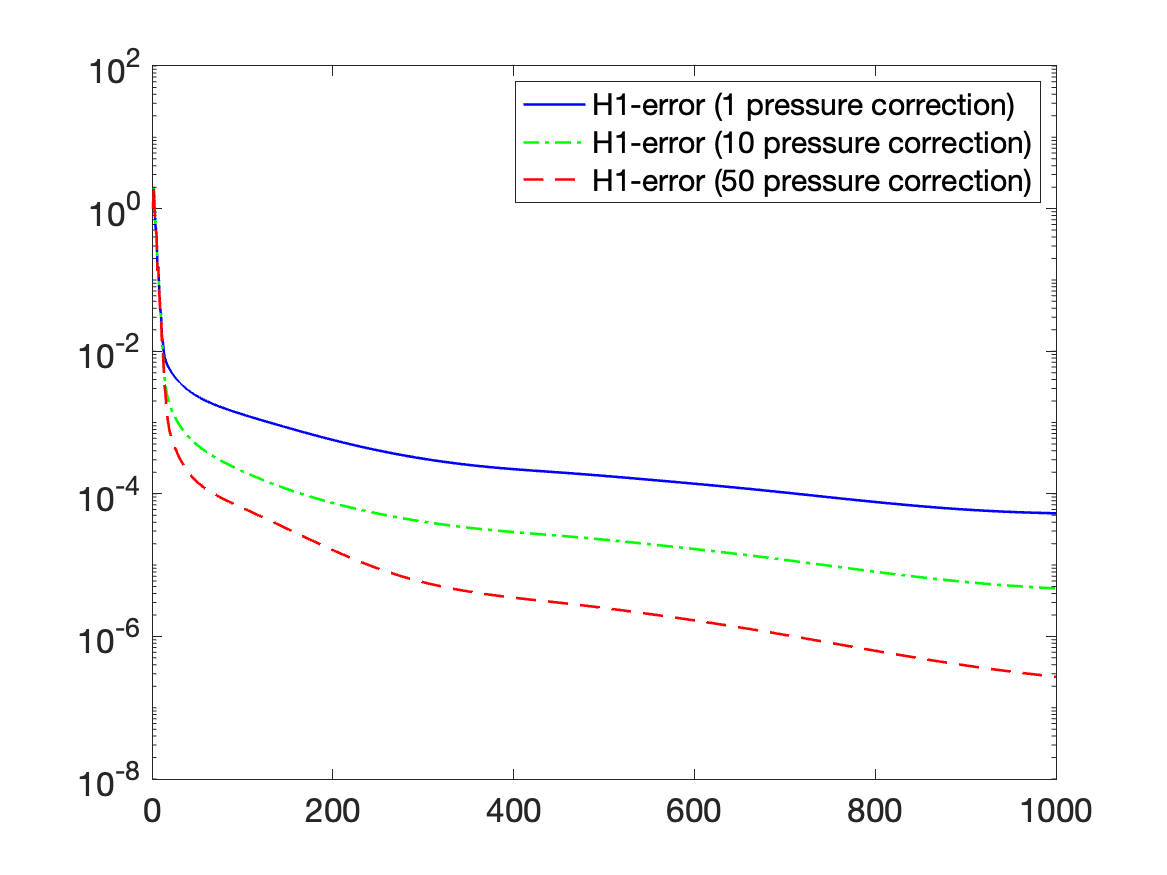}
\caption{$H^1$ rel. est. error (log-scale) over time (1000 timesteps, $\Delta t=0.01$): output and observer generated by FEM-solver.}
\label{fig:FEM}
\end{figure}
\begin{figure}
  \centering
  \includegraphics[width=1.0\columnwidth]{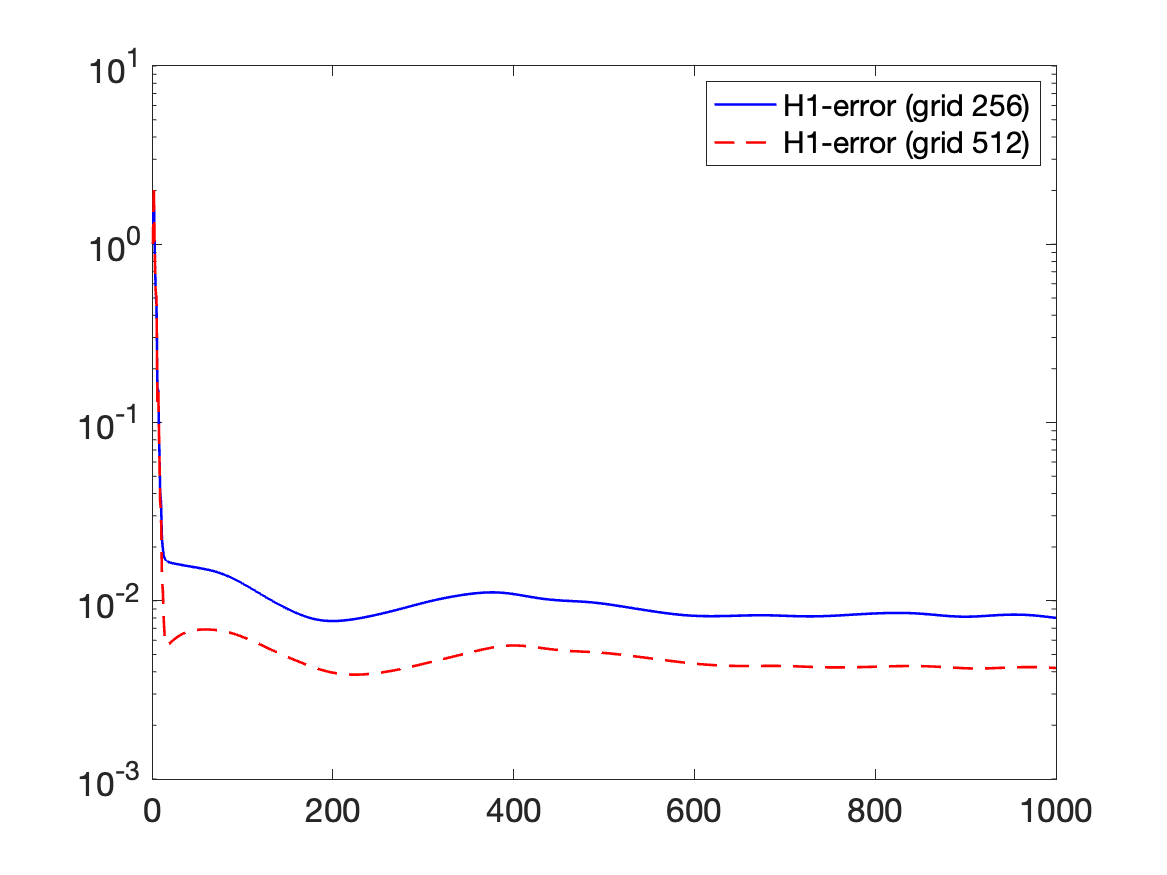}
  \caption{$H^1$ rel. est. error (log-scale) over time (1000 timesteps, $\Delta t=0.01$): output generated by FFT-solver, observer -- by FEM-solver.}
  \label{fig:observer}
\end{figure}

In what follows we first illustrate sensitivity of NSE with destabilizing Kolmogorov input to small perturbations of initial conditions. Then we illustrate convergence of the observer for the case of known inputs. And finally, we perform a crash-test: we take observations generated by one numerical method, and use them in the observer discretised by a different method.\\
For the crash-test we generate spatial averages outputs by a numerical solver, referred to as FFT-solver. FFT-solver is a pseudo-spectral numerical method, which relies upon vorticity-streamfunction formulation of NSE. It is exactly divergence free (as required by continuous formulation) and has spectral convergence property in space: its convergence rate in $H^1$ automatically increases with the degree of smoothness of initial conditions and inputs. For time discretisation we used 2nd order implicit midpoint with 5 iterations; an open-source implementation of FFT-solver with different time-stepping is available in \textit{jax-cfd} package\footnote{https://github.com/google/jax-cfd}. Then, we discretise the observer by a less accurate solver, referred to as FEM-solver. FEM-solver is implemented using Finite Element Method (based on Oasis Python package\footnote{https://github.com/mikaem/Oasis}~\cite{Mortensen15Oasis}) with 2nd order triangular elements providing global 1st order convergence rate in space. It also employs 2nd order Backward Differencing scheme for time discretisation. FEM-solver is not divergence free as it relies upon iterative minimization of velocity divergence at every time step. The immense differences between those solvers are pronounced on finite grids used below and their impact on observer, discretized by FEM-solver, is described as an unknown bounded disturbance $\vec d$.\\
        \textbf{Experiment setup.} We take a shifted domain $\Omega = [0,\ell_1]\times[0,\ell_2]$ with $\ell_1=\ell_2=2\pi$, and a destabilizing input is taken to be $\vec g(x,y) = [-5\sin(10y),0]^\top$. The initial velocity $\vu(0)$ is generated randomly and taken such that $\|\vu(0)\|_{H^1}\approx 9.6$. Both solvers do $1000$ steps forward in time with timestep $\Delta t=0.01$. Spatial resolution varies as detailed below.\\
        \textbf{Turbulent behaviour.} Top panel of Fig.~\ref{fig:TurbulenceA} illustrates the sensitivity of NSE to small perturbations of the Initial Condition (IC) measured in $H^1$-norm: red curve shows dynamics of $H^1$-distance between two trajectories, $\vu$ and $\vec v$ obtained by high-precision FFT-solver on $256\times256$-grid with $\nu=0.01$ and the same input $\vg$; $\vu$ and $\vec v$ are close by initially, $\|\vu(0)-\vec v(0)\|_{H^1}<10^{-4}\|\vu(0)\|_{H^1}$. If we repeat the same simulation but for $\nu=0.1$ NSE becomes stable: $H^1$-distance between two trajectories decays (blue curve). Bottom panel of Fig.~\ref{fig:TurbulenceA} shows dynamics of $H^1$-norm of two trajectories with same ICs and input but different $\nu$: for $\nu=0.1$ $H^1$-norm levels off, and the flow is laminar (stable) as shown on Fig.~\ref{fig:nu1}, in contrast, for $\nu=0.01$ $H^1$-norm is changing and the flow is turbulent as shown on Fig.~\ref{fig:nu01}.\\
\textbf{Known input.} In this test we show that detectability conditions of Theorem~\ref{t:1} are indeed simulation friendly: $\oC$ is taken to be spatial averages over squares $\Omega_j$ covering $\Omega$, and $h^2$ in~\eqref{eq:CPoin} is the area of the largest $\Omega_j$, thus $h^2$ in fact determines the number of squares (sensors) required for convergence. $h^2$ is found from $C1$ (Theorem~\ref{t:1}): since $\vec d=0$ it follows that we can set $R=0$, $\beta=1$ (as per Remark~\ref{sec:observer-design-0}). Also $\Sigma_T(s)=0$ for any $T>0$. Hence~\eqref{eq:tstarTeps} holds for every $T_1=T>0$ and $\varepsilon>0$ such that $2/(\nu\lambda_1) + 2\nu e^{-\lambda_1 \nu t^\star} \|\nabla\vu_0\|^2_{\vsLt}C_{\vg}^{-2} \le \varepsilon T_1 $ as $\delta(t^\star)=0$. We set $\varepsilon=0.001$ and $T_1= 1.01\varepsilon \times 2/(\nu\lambda_1)$ so~\eqref{eq:tstarTeps} holds for $t^\star$ such that 2nd term of the last inequality is less than $0.01\times 2/(\nu\lambda_1)$. We plug $\kappa=1+\varepsilon$ into~\eqref{eq:Theta} and maximize $\Theta$ by using grid search: we find $\Gamma_{\max}=611$. Hence $h^2=0.0029$ as per $C1$, and $L=262$. To get outputs we discretize NSE by FEM-solver with quadratic triangular elements constructed on a uniform grid of $256\times256$ nodes. Then the outputs are plugged into the observer discretized by the same FEM-solver. The $H^1$-estimation error is given in Fig.~\ref{fig:FEM} for 1,10 and 50 pressure corrections, which are used in FEM-solver to minimize the numerical divergence: 50 corrections have smallest divergence (red curve). Clearly, reduction of the numerical divergence implies reduction of $H^1$-estimation error.\\
\hspace{-73pt}
\textbf{Uncertain input.} In this test we pick $\varepsilon$, $T_1$ and $h,L$ as above but use FFT-solver to generate the outputs. The differences between discrete NSEs obtained by FEM-solver and FFT-solver on finite grids are significant, and in fact the former can be seen as the latter but with an additive uncertain input $\vec d$ which is expected to get ``smaller'' for finer grids. And this is exactly what we see on Fig.~\ref{fig:observer}: due to the presence of the disturbance the observer converges into a ``zone'' which shrinks (in $H^1$-norm) when spatial resolution increases from $256\times256$ to $512\times512$.
\hspace{-50pt}

\section{Conclusions}
\label{sec:conclusions}

We proposed {\it simulation friendly} detectability conditions for 2D Navier-Stokes Equation, and designed infinite-dimensional globally converging Luenberger observer for continuous in time measurements. Promising research directions include extending our approach to sampled measurements, and to pointlike measurement where $\Omega_j$ do not necessarily cover the entire domain $\Omega$ as studied in~\cite{Selivanov} for the case of 2D heat equation.




\bibliographystyle{plain}
\bibliography{refs/refs,refs/myrefs,refs/cdc_refs,rev1refs}

\appendix

\section{Proofs}
\label{sec:proofs}
\def\Elproofname{Proof of Lemma~\ref{sec:upper-bounds-linfty}}
\begin{proof}
Take $\vu=\bigl[
\begin{smallmatrix}u & v
\end{smallmatrix}
\bigr]\in [H_p^2(\Omega)]^2$. By Sobolev embedding theorems (see~\cite[Th.11.3.14]{Hutson1980}) $H^2(\Omega)$ is continuously embedded into the space of continuous functions $C(\Omega)$ hence the vector-function $\vu$ is continuous. Since $\vu$ is also periodic by Weierstrass approximation theorem~\cite[Corollary 3.1.11.]{Grafakos2008} $\vu$ can be represented as a uniform limit of its Fourier series:
\begin{equation}
  \label{eq:u-FS}
\lim_{N\to\infty}\sup_{x\in\Omega}\|\vu(\vx)-\sum_{\|\vk\|_{\R^2} \le N}\bigl[
\begin{smallmatrix}
  u_{\vk}&v_{\vk}
\end{smallmatrix}
\bigr] \ee^{2\pi\ii \frac{\vk}{\vL}\cdot \vx}\|_{\R^2}= 0
\end{equation}
with $u_{\vk}$ and $v_{\vk}$ being complex conjugates of $u_{-\vk}$ and $v_{-\vk}$ (since $\vu$ is real-valued), and
$u_{\vk}=(\ll_1\ll_2)^{-1}\int_\Omega e_1\cdot \vu(\vx)\ee^{-2\pi\ii \frac{\vk}{\vL}\cdot \vx} d\vx$,
and $v_{\vk}$ defined as $u_{\vk}$ but with $e_2$, provided $e_{1,2}$ -- canonical basis of $\R^2$. We have:
\begin{align}
  \|\vu(\vx)\|_{\R^2}&\le \|\vu(\vx)-\sum_{\|\vk\|_{\R^2} \le N}\bigl[
\begin{smallmatrix}
  u_{\vk}&v_{\vk}
\end{smallmatrix}
\bigr] \ee^{2\pi\ii \frac{\vk}{\vL}\cdot \vx}\|_{\R^2}     \label{eq:uxupb}\\
&+\sum_{\|\vk\|_{\R^2}\le N} (|u_{\vk}|_{\C}^2 + |v_{\vk}|_{\C}^2)^{\frac12}
|\ee^{2\pi\ii \frac{\vk}{\vL}\cdot \vx}|_{\C}\\
                   &\le \sum_{\vk\in\Z^2} (|u_{\vk}|_{\C}^2 + |v_{\vk}|_{\C}^2)^{\frac12}\label{eq:ul1}
\end{align}
To get~\eqref{eq:ul1} one sends $N\to\infty$, invokes~\eqref{eq:u-FS} and recalls that $|\ee^{2\pi\ii \frac{\vk}{\vL}\cdot \vx}|_{\C}=1$ and that $\vu\in L^2(\Omega)^2$, as $\vu$ is continuous, hence the series in~\eqref{eq:ul1} is converging.\\Let us compute 
$ \|\vu\|^2_{L^2} = \|u\|^2_{L^2}+ \|v\|^2_{L^2}$, $ \|\nabla \vu\|^2_{L^2} = \|\nabla u\|^2_{L^2}+ \|\nabla v\|^2_{L^2}$, $ \|\Delta \vu\|^2_{L^2} = \|\Delta u\|^2_{L^2}+ \|\Delta v\|^2_{L^2}$. Recall Parseval's identity
\begin{equation}
  \|\vu\|^2_{L^2}= \ll_1\ll_2 \sum_{\vk\in\Z^2} (|u_{k_1}|^2_\C +|u_{k_2}|^2_\C), \label{eq:u-l2}
\end{equation}
and classical relations between the smoothness of a function and decay of its Fourier coefficients (see~\cite[Theorem 3.2.9]{Grafakos2008}), and differentiate Fourier series of a function of $H^2$-class (e.g.~\cite[p.182]{Grafakos2008}) to compute norms of $\nabla\vu$ and $\Delta\vu$:

\[
  \|\nabla \vu\|_{L^2}^2 = 4\pi^2\ll_1\ll_2\sum_{\|\vk\|_{\R^2}\in\Z^2} \|\frac{\vk}{\vL}\|^2_{\R^2}(|u_{\vk}|_{\C}^2 + |v_{\vk}|_{\C}^2)
\]
\begin{equation}
   \|\Delta \vu\|_{L^2}^2 = 16\pi^4\ll_1\ll_2\sum_{\|\vk\|_{\R^2}\in\Z^2} \|\frac{\vk}{\vL}\|^4_{\R^2}(|u_{\vk}|_{\C}^2 + |v_{\vk}|_{\C}^2). \label{eq:u-Lap}
\end{equation}
 Now we split~\eqref{eq:ul1} into a finite sum and the remainder:
\begin{align*}
  \|\vu(\vx)\|_{\R^2}&\le \sum_{\|\vk\|_{\R^2}\le \gamma} \frac{\bigl(\ll_1\ll_2 (1+ 4\pi^2 \|\frac{\vk}{\vL}\|^2_{\R^2}) (|u_{\vk}|_{\C}^2 + |v_{\vk}|_{\C}^2)\bigr)^{\frac12}}{(\ll_1\ll_2 + 4\pi^2\ll_1\ll_2 \|\frac{\vk}{\vL}\|^2_{\R^2})^\frac12}\\
   &+\sum_{\|\vk\|_{\R^2}> \gamma}\frac{4\pi^2\sqrt{\ll_1\ll_2}\|\frac{\vk}{\vL}\|^2_{\R^2}}{4\pi^2\sqrt{\ll_1\ll_2}\|\frac{\vk}{\vL}\|^2_{\R^2}} (|u_{\vk}|_{\C}^2 + |v_{\vk}|_{\C}^2)^{\frac12}\\
\end{align*}
and make use of~\eqref{eq:u-l2}-\eqref{eq:u-Lap} to derive~\eqref{eq:AgmoN}:
\begin{align*}
   \|\vu(\vx)\|_{\R^2}&\le \|\vu\|_{H^1}\left(\sum_{\|\vk\|_{\R^2}\le \gamma}(\ll_1\ll_2 + 4\pi^2\ll_1\ll_2 \|\frac{\vk}{\vL}\|^2_{\R^2})^{-1}\right)^{\frac12}\\
   &+\|\Delta\vu\|_{L^2} \left(\sum_{\|\vk\|_{\R^2}> \gamma}(16\pi^4\ll_1\ll_2\|\frac{\vk}{\vL}\|^4_{\R^2})^{-1}\right)^{\frac12}\\
   &\le \|\vu\|_{H^1} \left(\int_{x_1^2+x_2^2\le \gamma} \frac{dx_1dx_2}{\ll_1\ll_2 + 4\pi^2\ll_1\ll_2 \|\frac{\vx}{\vL}\|^2_{\R^2}}\right)^{\frac12}\\
   &+\|\Delta\vu\|_{L^2} \left(\int_{x_1^2+x_2^2> \gamma} \frac{dx_1dx_2}{16\pi^4\ll_1\ll_2\|\frac{\vx}{\vL}\|^4_{\R^2}}\right)^{\frac12}\\
  &=\|\vu\|_{H^1}  I_1^{\frac12} + \|\Delta\vu\|_{L^2} I_2^{\frac12}.
\end{align*}
By coarea formula:
\[I_1=\int_0^\gamma dr \int_{x_1^2 + x_2^2 = r^2 }\frac{dS}{\ll_1\ll_2 + 4\pi^2\ll_1\ll_2 \|\frac{\vx}{\vL}\|^2_{\R^2}}\]
To compute the latter line intergal along the circle set $x_1=r\sin(t)$ and $x_2=r\cos(t)$:
\begin{align*}
  I_1 &=\int_0^\gamma \int_0^{2\pi} \frac{dr dt}{\ll_1\ll_2 + \frac{4\pi^2}{\ll_1\ll_2}(\ll_2^2 r^2 \sin^2(t)+\ll_1^2 r^2 \cos^2(t))}\\
      &=\int_0^\gamma \frac{dr}{(\ll_1^2+4\pi^2r^2)^{\frac12}(\ll_2^2+4\pi^2r^2)^{\frac12}}\\
      &\le \int_0^\gamma \frac{2\pi rdr}{4\pi^2 r^2+\ll_1\ll_2}=(4\pi)^{-1} \log(1+\frac{4\gamma^2\pi^2}{\ll_1\ll_2})
\end{align*}
Analogously we compute $I_2=\frac{\ll_1^2+\ll_2^2}{32\gamma^2 \pi^3}$. 
\end{proof}
\def\Elproofname{Proof of Prop.~\ref{prop:detectability-nse-1}}
\begin{proof}
  Take $\oC\in\DtctCls$ and set $\ve=\vu-\vz$. Subtracting \eqref{eq:NSEz} from~\eqref{eq:NSEu} we get ``the error equation'':  \[
  \dfrac {d}{dt}(\vec e,\vphi) + b(\vec e,\vec u,\vphi) + b(\vec z,\vec e,\vphi) + \nu((\vec e,\vphi)) = (\vf-\vF,\vphi)\eqno(*)
\]
Note that by~\eqref{eq:Poincare} $\|\nabla\ve\|^2_{\vsLt}=((\ve,\ve))\to0$ implies that $\|\ve\|^2_{\vsLt}=(\ve,\ve)\to 0$ hence, it is sufficient to demonstrate that
$\frac1T\int_t^{t+T} \|\oC\ve(\cdot,s)\|_{\vsLt}ds\to0$ and $\frac1T\int_t^{t+T} \|\vec f(\cdot,s)-\vec F(\cdot,s)\|_{\vsLt}^2 ds \to0$ (conditions $A)$ and $B)$ of Def.~\ref{def:detectability}) imply $V=((\ve,\ve))\to0$ if~\eqref{eq:h-bnd-detect} holds. To this end we plug $\vphi=A\vec e$ into the ``error equation'' $(*)$ and apply simple transformations: (i) recalling from~\eqref{eq:Poincare} and~\eqref{eq:gradA-bnd} that $(\vec e,A\vec e)=((\vec e,\vec e))$, $((\vec e,A\vec e))= (A\vec e,A\vec e)$, and (ii) recalling from~\eqref{eq:b:ort} that $b(\vec e,\vec e, A\vec e)=0$ which implies \(
b(\vec z,\vec e, A\vec e)=b(\vec u-\vec e,\vec e, A\vec e) = b(\vec u,\vec e, A\vec e)
\) so that
\begin{align*}
  b(\vec e,\vec u, A\vec e) &+ b(\vec z,\vec e, A\vec e)
  = b(\vec e,\vec u, A\vec e) + b(\vec u,\vec e, A\vec e) \\
  &= -b(\vec e,\vec e, A\vec u) \text{ (see \cite[F.(A.63)]{NSE-Turbulence-book2001})}
\end{align*}
we get:
\begin{equation}
      \label{eq:error-Ae}
      \dfrac{d}{dt}((\vec e,\vec e)) + \nu(A\vec e,A\vec e) = (\vf-\vec F,A\vec e) + b(\vec e,\vec e, A\vec u)
\end{equation}
Let us demonstrate that $\nu(A\vec e,A\vec e)$ ``dominates'' $ b(\vec e,\vec e, A\vec u)$ provided that condition $A)$ holds. Indeed
\begin{align*}
     b(\vec e,\vec e, A\vec u) &\overset{\eqref{eq:AgmoN}}{\le} C_1  \log^{\frac12}\bigl(1+\frac{4\pi^2 \gamma^2}{\ll_1\ll_2}\bigr)\|\nabla\ve\|^2_{\vsLt}\|A\vu\|_{\vsLt} \\
                               &+ C_2 \lambda_1^{-1}\gamma^{-1} \|A \vu\|_{\vsLt} \|A\ve\|^2_{\vsLt}
\end{align*}
Now, by~\eqref{eq:InterpIneq}, \eqref{eq:CPoin} and~\eqref{eq:gradA-bnd} we get:
\begin{align}
      \|\nabla\ve\|^2_{\vsLt}
      &\le  C_{\nabla} (\|\oC\ve\|_{\vsLt}\|A\ve\|_{\vsLt} + h C^{\frac12}_\Omega \lambda_1^{-\frac12} \|A\ve\|^2_{\vsLt})\label{eq:nabla-interp-bnd}
\end{align}
Set $\gamma=\|A \vu\|_{\vsLt} \Gamma$ for some $\Gamma>0$ and define
\begin{align*}
	\tilde C^1_{h,\Gamma}(t) &= C_1  C_{\nabla} \log^{\frac12}\bigl(1+\frac{4\pi^2 \|A \vu\|_{\vsLt}^2 \Gamma^2}{\ll_1\ll_2}\bigr) h C^{\frac12}_\Omega \lambda_1^{-\frac12}\|A \vu\|_{\vsLt}\\
	&+ C_2 \lambda_1^{-1}\Gamma^{-1}\\
	\tilde C_\Gamma(t) &= C_1  C_{\nabla} \log^{\frac12}\bigl(1+\frac{4\pi^2 \|A \vu\|_{\vsLt}^2 \Gamma^2}{\ll_1\ll_2}\bigr)
\end{align*}
Using $\tilde C^1_{h,\Gamma}$ and $\tilde C_\Gamma$ and noting that $2\|A\ve\|_{\vsLt}\|A\vu\|_{\vsLt}\le \|A\ve\|^2_{\vsLt}+\|A\vu\|^2_{\vsLt}$ we transform the upper bound for $b$:
\begin{equation}\label{eq:bUp}
	\begin{split}
		b(\vec e,\vec e, A\vec u)&\le \tilde C^1_{h,\Gamma}(t) \|A\ve\|^2_{\vsLt}\\
        &+ 0.5\tilde C_\Gamma(t) \|\oC\ve\|_{\vsLt}(\|A\ve\|^2_{\vsLt}+\|A\vu\|^2_{\vsLt})
  \end{split}
\end{equation}
By Schwartz inequality:
$$(\vf-\vF,A\ve)\le 1/\nu\|\vf-\vF\|_{\vsLt}^2 + \nu/4(A\ve,A\ve)$$
We plug the latter inequality and~\eqref{eq:bUp} into~\eqref{eq:error-Ae} (recall that  $V=((e,e))$):
\begin{align*}
                  &\dot V + (3\nu/4- \tilde C^1_{h,\Gamma}(t)) (A\vec e,A\vec e)\le\beta(t)\\
                  &\beta:=\|\vf-\vF\|_{\vsLt}^2/\nu
                  + 0.5\tilde C_\Gamma(t)\|\oC\ve\|_{\vsLt}(\|A\ve\|^2_{\vsLt}+\|A\vu\|^2_{\vsLt})
                \end{align*}
                If we set $\alpha(t)=(3\nu/4- \tilde C^1_{h,\Gamma}(t))$ then $V$ verifies the inequality $\dot V(t) + \alpha(t) V(t) \le \beta(t)$. To show that $V\to0$ we employ Lemma~1.1 from~\cite[p.125]{NSE-Turbulence-book2001}, a generalisation of the classical Gronwall lemma which in our case reads as follows: if $V$ verifies $\dot V(t) + \alpha(t) V(t) \le \beta(t)$ with the just defined $\alpha,\beta$ then $\lim_{t\to\infty}V(t) = 0$ provided there exist $T>0$ such that $\limsup_{t\to\infty}\frac1T\int_t^{t+T} \beta(s) ds =0$ and $\liminf_{t\to\infty}{\frac1T\int_t^{t+T} \alpha(s) ds} >0$. We claim that conditions $A)$ and $B)$ of Def.~\ref{def:detectability} imply the aforementioned condition for $\beta$, and~\eqref{eq:h-bnd-detect} imply the required condition on $\alpha$. To show the former recall from~\cite[p.101, A.60]{NSE-Turbulence-book2001} that $\sup_{t>t_\Delta}\|A\vu(t)\|_{\vsLt}<C_\Delta$ for some $C_\Delta >0$ which depends on $\|f\|_{\vsLiHza}$. Similar bound holds for $A\vz$ but depends on $\|F\|_{\vsLiHza}$. Hence $A\ve = A\vu-A\vz$ is bounded for $t>t_\Delta$ and so is $\tilde C_\Gamma(t)$. If conditions $A)$ and $B)$ of Def.~\ref{def:detectability} hold then $\frac1T\int_t^{t+T}\beta(s) ds\to 0$ as $t\to\infty$. Now, to demonstrate condition on $\alpha$ let us show that~\eqref{eq:h-bnd-detect} implies
                \begin{equation}
                  \label{eq:C1hnu}
\frac1T\int_t^{t+T} \tilde C^1_{h,\Gamma}(s) ds <(3\nu)/4
\end{equation}
Indeed, as it follows from~\eqref{eq:Au-ThetatT} and \eqref{eq:Au-L2avrg-bnd} for any $1<\kappa<<2$ there exist $t^\star>0$ and $T>0$ such that \[
  \frac1T\int_t^{T+t}\|A\vu\|^2_{\vsLt}ds  \le \theta_{t,T}\le  \kappa\|\vf\|^2_{\vsLiHza}\nu^{-2}, t>t^\star
\]
Having this in mind and applying integral Young inequality for
concave function
\begin{equation}
  \label{eq:Young}
  \frac1T\int_t^{t+T}\log(\star)ds\le \log(\frac1T\int_t^{t+T}(\star)ds)
\end{equation}
we compute:
\begin{align*}
  \frac1T&\int_t^{t+T} \tilde C^1_{h,\Gamma}(s) ds \le C_2 \lambda_1^{-1}\Gamma^{-1} + C_1 C_\nabla h  C^{\frac12}_\Omega \lambda_1^{-\frac12} \theta_{\nu,T}^{\frac12} \times\\
  &\quad \times\bigl(\frac1T\int_t^{t+T}\log\bigl(1+\frac{4\pi^2 \Gamma^2 \|A \vu(t)\|_{\vsLt}^2}{\ll_1\ll_2}\bigr)ds\bigr) ^{\frac12}\\
         &\le C_2 \lambda_1^{-1}\Gamma^{-1} + C_1 C_\nabla h  C^{\frac12}_\Omega \lambda_1^{-\frac12}\kappa^{\frac12}\|\vf\|_{\vsLiHza}/\nu \times\\
  &\qquad \times\log^{\frac12}\bigl(1+\frac{4\pi^2 \kappa \|\vf\|_{\vsLiHza}^2\Gamma^2 }{\nu^2\ll_1\ll_2}\bigr)
\end{align*}
This demonstrates that \eqref{eq:h-bnd-detect} implies~\eqref{eq:C1hnu}.\\
If $\oC\in\DtctClsA$ one needs to invoke~\eqref{eq:CPoin2} in~\eqref{eq:nabla-interp-bnd}, and~\eqref{eq:gradA-bnd} is not required so $\lambda_1^{-\frac12}$ in the r.h.s. of~\eqref{eq:nabla-interp-bnd} disappears as well as in the numerator of~\eqref{eq:h-bnd-detect}.
\end{proof}
\def\Elproofname{Proof of Lemma~\ref{l:observer-existence}}
\begin{proof}
  We employ the standard argument~\cite[p.23,§3.3]{TemamNSE-Periodic} with a difference in energy bounds due to the term $L\oC\ve$ which we outline below.
  Consider $\oC\in\DtctCls$. Galerkin projection of~(\ref{eq:NSE-LO}) onto $W_m$, the span of $m$ eigen-functions of the Stokes operator $A$, is obtained by substituting $\vz$ with $z_m$, the projection of $\vz$ onto $W_m$, and restricting test-functions to $\vec\phi\in W_m$, specifically for $\vec \phi=Az_m$ one gets: $\dfrac{d}{dt}((z_m,z_m)) + \nu\|Az_m\|^2_{\vsLt} = (\vg+L\oC\ve,Az_m)$. Adding and substructing $L(z_m,Az_m)$ and invoking~\eqref{eq:CPoin} after simple manipulations one finds:  
  \begin{equation}
    \label{eq:NSE-GalProj}
     \dfrac{d}{dt}((z_m,z_m)) + \nu/4(Az_m,Az_m) \le  1/\nu\|\vg+L\oC\vu\|^2_{\vsLt} 
   \end{equation}
   provided $L h^2 C_\Omega/2\le \nu$. Note that~\eqref{eq:NSE-GalProj} is similar to the classical a-priori energy bound for 2D NSE with periodic BC, e.g.~\cite[p.102,(A.65)]{NSE-Turbulence-book2001}. Then taking $m\to\infty$ and using the compactness argument~\cite[p.23,§3.3]{TemamNSE-Periodic} one deduces lemma's statement. The case of $\oC\in\DtctClsA$ follows the same logic.
\end{proof}
\def\Elproofname{Proof of Theorem~\ref{t:1}}
\begin{proof}
  Note that if $L$ and $h$ verify $C1$ or $C2$ then there is the unique $\vz\in L^\infty(\R_+,\vsVza)\cap L^2(t_0,t_1,\dom(A))$, $0\le t_0<t_1<+\infty$ (Lemma~\ref{l:observer-existence}). Recall that for $\vz$ and $\vu$ solving~\eqref{eq:NSEz} and~\eqref{eq:NSEu} with generic $\vF$ and $\vf$ respectively the dynamics of $((\ve,\ve))$, $\ve=\vu-\vz$ is governed by~\eqref{eq:error-Ae}. Now, let $\vu$ solve~\eqref{eq:NSEu} with $\vf=\vg+\vec d$, and $\vz$ solve~\eqref{eq:NSE-LO}, and plug $\vF = \vg+L\oC\ve$ and $\vf=\vg+\vec d$ into~\eqref{eq:error-Ae}: the resulting equation will determine the dynamics of $V=\|\nabla\ve\|^2_{\vsLt}=((\ve,\ve))$ defined in theorem's statement. With this in mind let us demonstrate that $C1)$ implies~\eqref{eq:ineq-V1}. To this end we transform~\eqref{eq:error-Ae}: for any $\Lambda_{1,2}>0$ and $q^2=2\|\vf-\vg\|^2_{\vsLt}=2\|\vec d\|^2_{\vsLt}$ %
   \begin{align}
     (\vf-\vec F,A\vec e) \le& (\vec d + L(\ve-\oC\ve),A\ve) - L(\ve,A\ve)\label{eq:vgLA}\\
                          &+ \Lambda_1(h^2 C_\Omega \|\nabla\ve\|_{\vsLt}^2 - \|\ve - \oC\ve\|_{\vsLt}^2) \label{eq:vgLA1}\\
                          &+ \Lambda_2(q^2(t) - \|\vec d\|_{\vsLt}^2) \label{eq:vgLA2}
   \end{align}
   where~\eqref{eq:vgLA1} is non-negative by~\eqref{eq:CPoin}. Recall definition of $C_{1,2}$ from~\eqref{eq:C1C2}, and~\eqref{eq:AgmoN}: for any $0<\beta<1$ we have
   \begin{align}
     b(\vec e,\vec e, A&\vec u) \le C_1  \log^{\frac12}\bigl(1+\frac{4\pi^2 \gamma^2}{\ll_1\ll_2}\bigr)\|\nabla\ve\|^2_{\vsLt}\|A\vu\|_{\vsLt}\label{eq:bAgmon} \\
                               &+ (1+\beta - \beta)C_2 \lambda_1^{-1}\gamma^{-1} \|A \vu\|_{\vsLt} \|A\ve\|^2_{\vsLt}\label{eq:bAgmon1}
\end{align}
Set $\gamma=\|A \vu\|_{\vsLt} \Gamma$ with $\Gamma>C_2/(\lambda_1\nu)$, define $\psi_{\Gamma,\nu}=\nu -C_2 /(\lambda_1\Gamma)$, add and substruct $\beta \nu \|A\ve\|^2_{\vsLt}$ to the l.h.s. of~\eqref{eq:error-Ae}, and substitute~\eqref{eq:vgLA}-\eqref{eq:bAgmon1} into r.h.s. of~\eqref{eq:error-Ae}. Collecting the terms with $A\ve$, $\ve-\oC\ve$ and $\vec d$, and transforming the resulting expressions to sums of squares after simple algebra we get:
\begin{align*}
  \dot V &\le \alpha(t)V + \Lambda^*_2q^2(t) -  \|(1-\beta)^\frac12 \psi_{\Gamma,\nu}^{\frac12} A\ve - \sqrt{\Lambda_2^*}\vec d \|^2_{\vsLt}\\
  &- \|\beta^{\frac12} \psi_{\Gamma,\nu}^{\frac12} A\ve - \sqrt{\Lambda_1^*}(\ve-\oC\ve)\|^2_{\vsLt}
\end{align*}
provided $\Lambda^*_1$ solves $\Lambda^*_1 4\beta \psi_{\Gamma,\nu}=L^2$, $\Lambda^*_2$ solves $\Lambda^*_2 4(1-\beta)\psi_{\Gamma,\nu}=1$ and
\begin{align}
  \alpha &= - L + \Lambda^*_1 h^2 C_\Omega+ \Psi(\vu,\Gamma)\label{eq:def:alpha}\\
  \Psi(\vu,\Gamma)&=C_1  \|A \vu\|_{\vsLt}\log^{\frac12}\bigl(1+\frac{4\pi^2 \Gamma^2 \|A \vu\|_{\vsLt}^2}{\ll_1\ll_2}\bigr) \label{eq:def:Psi}
\end{align}
Hence, $\dot V \le \alpha(t)V + \Lambda^*_2q^2(t)$ and by classical Gronwall lemma we have for $t\ge s$:
\begin{equation}
  \label{eq:GBl}
V(t)\le V(s) e^{\int_s^t \alpha(\tau)d\tau} + \Lambda^*_2\int_s^t  q^2(\tau)  e^{\int_\tau^t \alpha(\sigma)d\sigma} d\tau
\end{equation}
Let us shows that 1st and 2nd terms in r.h.s. of~\eqref{eq:GBl} are bounded by the 1st and 2nd terms of~\eqref{eq:ineq-V1} respectively. To bound $e^{\int_s^t \alpha d\tau}$ we show that if $t^\star$ and $T_1$ verify~\eqref{eq:tstarTeps} then for $W(L,\Gamma)= - L + L^2 h^2 /(2\beta \hat L_\nabla(\Gamma)) + \hat L_\Delta(\Gamma)$ with $\hat L_{\nabla,\Delta}$ defined in~\eqref{eq:hatL}-\eqref{eq:hatL1} it holds:
\begin{align}\label{eq:alphaneg}
\forall t>t^\star: \quad \frac1{T_1}\int_t^{t+T_1}\alpha(s)ds \le W(L,\Gamma)
\end{align}
Indeed, let us take $T_1\ge T$ such that the 2nd term of~\eqref{eq:tstarTeps} is small (e.g. less than $\varepsilon/2$), then note that~\eqref{eq:fgdecay} implies $\lim_{s\to\infty}\Sigma_{T_1}(s)=0$ for any $T_1\ge T$, hence taking a large enough $t^\star$ one can make $\delta(t^\star)$ (1st term of~\eqref{eq:tstarTeps}) and last term of~\eqref{eq:tstarTeps} small enough (e.g. less than $\varepsilon/2$) for~\eqref{eq:tstarTeps} to hold. For these $T_1$ and $t^\star$ we bound $\Psi(\vu,\Gamma)$ defined in~\eqref{eq:def:Psi}: to this end recall Young inequality~\eqref{eq:Young} which together with~\eqref{eq:Au-L2avrg-bnd} implies for any $T_1>0$
\begin{align}
  \frac1{T_1}\int_t^{t+T_1}\Psi(\vu,\Gamma)ds &\le
  C_1  \theta_{t,T_1}^{\frac12}\log^{\frac12}\bigl(1+\frac{4\pi^2 \theta_{t,T_1}\Gamma^2 }{\ll_1\ll_2}\bigr)
  \label{eq:theta:ubnd}
\end{align}
with $\theta_{t,T_1}$ defined in~\eqref{eq:Au-ThetatT}. Now, we bound $\theta_{t,T_1}$: recall definitions of $\kappa$ and $\ball$, and plug $\vf = \vg+\vec d$ into~\eqref{eq:Au-ThetatT}. Noting that $\|\vf\|_{\vsLiHza}\le R+\gnrm$ as $\vec d\in\ball$ we find that 1st term in~\eqref{eq:Au-ThetatT} is bounded by the 2nd term of~\eqref{eq:tstarTeps}, and, by Cauchy-Schwartz inequality, the 2nd term in~\eqref{eq:Au-ThetatT} is bounded by $\frac{\gnrm^2}{\nu^2}+\delta(t^\star)$ hence by~\eqref{eq:tstarTeps}:
    \begin{equation}
    \label{eq:Thetabnd}
    \begin{split}
      \theta_{t,T_1} &\le \bigl(\frac{2(\gnrm+R)^2}{\gnrm^2 T_1\nu\lambda_1} + \frac{2\nu\mathbb{e}^{-\lambda_1 \nu t}\|\nabla\vec u_0\|^2_{\vsLt}}{T_1 \gnrm^2} +\delta\bigr)\times\\
      &\times\frac{\gnrm^2}{\nu^2}+ \frac{\gnrm^2}{\nu^2}\le  \kappa \frac{\gnrm^2}{\nu^{2}},\,t>t^\star, T_1\ge T
    \end{split}
  \end{equation}
  Substituting~\eqref{eq:def:alpha}, \eqref{eq:theta:ubnd} and~\eqref{eq:Thetabnd} into l.h.s. of~\eqref{eq:alphaneg}, noting that $\Lambda^*_1=L^2/(4\beta\psi_{\Gamma,\nu})$ and recalling definition of $\hat L_\nabla$ we obtain~\eqref{eq:alphaneg}.

Let us show that $W(L,\Gamma)<0$ provided $h$ and $L$ are chosen as in $C1)$. Indeed, $W$ is a quadratic polynomial (in $L$) with two distinct real roots \(
  L_\pm=(1\pm\sqrt{1-4  a \hat L_\Delta}\,)/2a\),
  $a = h^2 /(2\beta \hat L_\nabla(\Gamma))$ provided the discriminant of $W$ is positive: $1-4a \hat L_\Delta> 0$. The latter implies $0<L_-<L_+$ and since $W(0,\Gamma) =\hat L_\Delta>0$ it follows that $W(L,\Gamma)<0$ for any $L\in(L_-,L_+)$. Hence, $W(L,\Gamma)<0$ for $L=\beta\hat L_\nabla /h^2= (L_++L_-)/2 \in (L_-,L_+)$ provided $4a \hat L_\Delta<1$ which is the case for $h^2< (2\beta \hat L_\nabla(\Gamma))/(4 \hat L_\Delta(\Gamma))=\beta \Theta(\Gamma)$ (compare to $C1)$). Specifically, $W(L, \Gamma_{\mathrm{max}} )=Q(L)<0$ for $Q(L)$ defined in~\eqref{eq:ineq-V1}, and $\omega$, $L$ chosen as in $C1)$, and in fact choosing $\Gamma = \Gamma_{\mathrm{max}}$ allows to maximize the upper-bound for $h^2$. Hence, $\forall t>t^\star$ and $\omega=1/2$, $ L =\beta\hat L_\nabla(\Gamma_{\mathrm{max}}) /h^2$:
\begin{align}
  \label{eq:alpha-QL}
\frac1{T_1} \int_t^{t+T_1}\alpha(s)ds  \le Q(L)<0
\end{align}
Let us now bound $e^{\int_s^t \alpha d\tau}$. To this end we follow~\cite[p.156, f.(A.4)]{NSE-Turbulence-book2001}, namely the argument given after formula (A.3) there: take $k=\lfloor\frac{t-s}{T_1}\rfloor$, that is the smallest integer $k\ge0$ such that $s+kT_1\le t\le s+(k+1)T_1$, and note that
\begin{equation}
  \label{eq:exp-alpha-bnd}
  \begin{split}
    &e^{\int_s^t \alpha(\tau)d\tau}=  e^{\int_s^{s+kT_1} \alpha(\tau)d\tau}e^{\int_{s+kT_1}^t \alpha(\tau)d\tau }\\
    &\overset{\eqref{eq:alpha-QL},\eqref{eq:def:alpha}}{\le} e^{Q(L)T_1 k } e^{-L\omega(t-s-kT_1)+\int_{s+kT_1}^{s+(k+1)T_1} \Psi(\vu,\Gamma_{\mathrm{max}})d\tau }\\
    &\le e^{Q(L) T_1(k+1)}e^{-T_1Q(L)}  e^{-L\omega(t-s-kT_1) + T_1 \hat L_\Delta(\Gamma_{\mathrm{max}})}\\
    &\le e^{(Q(L)-\omega L)(t-s) + \omega L T_1\lfloor\frac{t-s}{T_1}\rfloor}e^{\omega LT_1}
  \end{split}
\end{equation}
where to go from 2nd to 3rd line of~\eqref{eq:exp-alpha-bnd} we employed \eqref{eq:theta:ubnd},\eqref{eq:Thetabnd},\eqref{eq:hatL1}, and from 3rd to 4th -- $e^{T_1 (k+1) Q(L)}\le e^{Q(L)(t-s)}$ (as $Q(L)<0$ and $t-s\le(k+1)T_1$).\\ Now, let us bound the 2nd term of~\eqref{eq:GBl}. Again, we follow the derivation of~\cite[p.156, f.(A.5)]{NSE-Turbulence-book2001}, namely the argument given after formula (A.4) there, which is based on exponential sum formulas:
\begin{equation}
  \label{eq:qbnd}
\int_s^t  q^2(\tau)e^{\int_\tau^t \alpha(\sigma)d\sigma} d\tau \le 2 T_1 e^{T_1 (\omega L - Q(L))}\Sigma_{T_1}(s)
\end{equation}
Finally, recalling~\eqref{eq:hatL} and definition of $\Lambda_2^\star$ we deduce~\eqref{eq:ineq-V1} from \eqref{eq:GBl}, \eqref{eq:exp-alpha-bnd} and \eqref{eq:qbnd}.

Now, let us demonstrate that $C2)$ implies~\eqref{eq:ineq-V1}. If $\oC\in\DtctClsA$ then~\eqref{eq:vgLA1} changes to $\Lambda_1(h^2 C_\Omega \|A\ve\|_{\vsLt}^2 - \|\ve - \oC\ve\|_{\vsLt}^2)$ as it follows from~\eqref{eq:CPoin2}. Now, as above (see discussion right after~\eqref{eq:bAgmon1}), we substitute~\eqref{eq:vgLA}, modified~\eqref{eq:vgLA1} and~\eqref{eq:vgLA2}-\eqref{eq:bAgmon1} into r.h.s. of~\eqref{eq:error-Ae}, and collect terms with $A\ve$, $\ve-\oC\ve$ and $\vec d$:
\begin{align*}
  \dot V &\le \alpha'(t)V + \Lambda^\star_2q^2(t) -  \|(1-\beta)^\frac12 \psi_{\Gamma,\nu}^{\frac12} A\ve - \sqrt{\Lambda_2^*}\vec d \|^2_{\vsLt}\\
         &- \|(\beta\psi_{\Gamma,\nu}-\Lambda'_1 h^2 C_\Omega)^{\frac12}A\ve - \sqrt{\Lambda_1'}(\ve-\oC\ve)\|^2_{\vsLt}
\end{align*}
provided $\Lambda^*_2$ solves $\Lambda^*_2 4(1-\beta)\psi_{\Gamma,\nu}=1$ (as above), and
\begin{equation}
  \label{eq:alpha-prime}
\alpha' =  - L + \Psi(\vu,\Gamma) 
\end{equation}
and $\Lambda_1'$ solves a quadratic equation $4\Lambda_1(\beta\psi_{\Gamma,\nu}-\Lambda_1 h^2 C_\Omega) = L^2$, solution of which is given by:
\begin{align}
  \Lambda_1' &= \frac{\psi_{\Gamma,\nu}+\sqrt{\beta^2\psi_{\Gamma,\nu}^2-L^2 C_\Omega h^2}}{2h^2 C_\Omega}\\
  &\text{provided: }\psi_{\Gamma,\nu}>0, \quad L h\sqrt{C_\Omega} < \beta\psi_{\Gamma,\nu}\label{eq:Lambdaprime}
\end{align}
Similarly, to demonstrate~\eqref{eq:ineq-V1} one can repeat all the steps given after \eqref{eq:def:Psi}: employ GB-lemma to derive analog of~\eqref{eq:GBl}, and then bound its terms. Since the last term is the same, the only bound that remains to be established is the analog of~\eqref{eq:alpha-QL}, namely that $\frac1{T_1}\int_t^{t+T_1}\alpha'(s)ds < Q(L)=-L\omega + \hat L_\Delta(\Gamma_{\mathrm{max}})<0$ for $L$ and $\omega$ defined in $C2)$. But the first part of the latter inequality immediately follows from~\eqref{eq:alpha-prime} and \eqref{eq:theta:ubnd}, \eqref{eq:Thetabnd}, \eqref{eq:hatL1}, and the inequality $Q(L)<0$ holds if $L=\theta\hat L_\Delta(\Gamma_{\mathrm{max}})$ for $\theta>1$ (as suggested in $C2)$). To satisfy~\eqref{eq:Lambdaprime} we set $\Gamma=\Gamma_{\mathrm{max}}$ and substitute $L=\theta\hat L_\Delta(\Gamma_{\mathrm{max}})$ into~\eqref{eq:Lambdaprime} and get a condition on $h: h < \frac{\beta\psi_{\Gamma,\nu}}{\theta \hat L_\Delta\sqrt{C_\Omega}}=\beta \sqrt{C_\Omega} \theta^{-1} \Theta(\Gamma)$ as suggested in $C2)$. Note that $\Gamma_{\mathrm{max}}>C_2/(\lambda_1\nu)$ hence $\psi_{\Gamma_{\mathrm{max}},\nu}>0$ as required by~\eqref{eq:Lambdaprime}. Hence,~\eqref{eq:ineq-V1} holds for $V$ if $L$, $\omega$ and $h$ are defined as in $C2)$. 
\end{proof}


\end{document}